\documentclass[journal,twoside,web]{ieeecolor}
\usepackage{generic}
\usepackage{cite}
\usepackage{amsmath,amssymb,amsfonts}
\usepackage{algorithmic}
\usepackage{graphicx}
\usepackage{algorithm,algorithmic}

\newtheorem{assumption}{Assumption}
\newtheorem{theorem}{Theorem}

\newtheorem{remark}{Remark}

\newtheorem{definition}{Definition}
\allowdisplaybreaks
\usepackage{multirow}
\usepackage{subcaption}
\usepackage{textcomp}
\def\BibTeX{{\rm B\kern-.05em{\sc i\kern-.025em b}\kern-.08em
    T\kern-.1667em\lower.7ex\hbox{E}\kern-.125emX}}
\begin{document}


\title{Asymptotic, Exponential, and Prescribed-Time Unbiasing in Seeking of Time-Varying Extrema}

\author{Cemal Tugrul Yilmaz, Mamadou Diagne, and Miroslav Krstic, \IEEEmembership{Fellow, IEEE}
\thanks{The authors are with the Department of Mechanical and Aerospace
Engineering, University of California San Diego, La Jolla, CA 92093-
0411 USA (e-mail: cyilmaz@ucsd.edu; mdiagne@ucsd.edu, krstic@ucsd.edu).}}

\maketitle

\allowdisplaybreaks

\begin{abstract}

Our recently developed ``unbiased'' extremum seeking (uES) algorithms ensure perfect convergence to the optimum at a user-assigned exponential rate, in spite of not using the Newton approach, or, more powerfully, within a user-prescribed time. Unlike classical approach, these algorithms use time-varying adaptation and controller gains, along with constant or time-varying probing frequencies (chirp signals). This paper advances our earlier uES designs from strongly convex maps with static optima to a broader class of convex cost functions with time-varying optima diverging at arbitrary rates, even in finite time.
This advancement first motivates the use of Lie bracket averaging instead of classical averaging due to the average system system, which doesn’t necessarily need to be exponentially convergent, and the existence of non-periodic time-varying parameters; second, it necessitates the formulation of non-trivial and key feasibility conditions for the choice of time-varying design parameters and their decay/growth rates in relation to the convexity of the map and the divergence rate of optima. These conditions indicate that, for constant-frequency probing, the user-defined asymptotic rate of unbiasing is limited by the convexity of the map. However, this rate can be made arbitrarily fast (including asymptotic, exponential, and prescribed time) using chirpy probing, which requires sufficiently rapid frequency and adaptation growth to enable tracking of faster-diverging optima.
Furthermore, our uES methodology leverages state and time scaling techniques. State scaling involves scaling the state by an unboundedly growing function, while time scaling involves time dilation/contraction transformations. The practical stability of the state-scaled system in the dilated time domain is established through Lie bracket averaging, which provides unbiased and accelerated convergence of states to the optimum in the original time domain. 
In addition to numerical simulations of the designs, we experimentally test the feasibility of exponential uES for tuning the angular velocity of a unicycle to seek a static light source.

\end{abstract}

\begin{IEEEkeywords}
unbiased extremum seeking, prescribed-time control, time-varying optima, unbiased source seeking
\end{IEEEkeywords}

\section{Introduction}
\label{sec:introduction}
``Unbiasing'' of extremum seekers, namely, perfect convergence to the extremum by a gradual reduction of the perturbation amplitude to zero, has recently emerged \cite{yilmaz2023press} but it has been limited to the basic idea of using an {\em exponential} rate of unbiasing for strongly convex (locally {\em quadratic}) maps. In a major effort in redesign and analysis, to enable unbiasing for maps that are strictly but {\em not strongly convex}, 
we provide algorithms and proofs for unbiasing that identify sufficiently slow rates of decay of the amplitude in relation to the ``flatness'' of the map at the extremum. This allows oscillations to fade slowly enough to ``explore'' around the optimum. ``Exploiting'' the gradient estimate with growing controller gain and the use of chirp probing opens the door to further acceleration in the rate of unbiasing, up to a user-prescribed time.

Practical motivation for unbiased seeking of extrema of maps that are not strongly convex abounds. Reducing or eliminating the residual oscillation enables increased power production in photovoltaic systems \cite{moura2013lyapunov}, more precise tuning of design parameters of feedback control laws \cite{wu2021active}, and more accurate localization of leakage sources by mobile robots \cite{jabeen2023robot}.


The importance of unbiasing extends beyond maps that lack strong convexity to maps whose extrema vary with time. Time-varying extrema occur in renewable energy systems \cite{krstic2014extremum,moura2013lyapunov}, where controller parameters need to be adjusted to maximize power generation despite fluctuating weather conditions and varying energy demand. Similarly, in mineral processing systems, continuously adjusting aeration rates is necessary to maintain optimal air recovery \cite{wepener2023extremum}. We show how persistently-changing extrema can be sought in a bias-free manner, provided the acceleration of the extremum decays at a sufficient rate.

\subsection{ES with unbiased convergence to fixed optima}
To address persistent oscillation around the fixed optimum, \cite{tan2009global} introduces a method involving a time-varying perturbation amplitude that decays to zero, ensuring practical asymptotic convergence to the global extremum despite local extrema. Another approach in \cite{wang2016stability} adjusts the perturbation amplitude based on system output. However, the claim of exponential convergence to the optimum in \cite{wang2016stability} has been theoretically \cite{atta2019comment} and numerically \cite{guay2020uncertainty, haring2016asymptotic} disputed. Kalman filtering techniques, with feedback to dynamically update the amplitude, are employed in \cite{pokhrel2023control} and \cite{bhattacharjee2021extremum}, achieving asymptotic convergence to a neighborhood of the extremum with diminishing oscillation.
In \cite{lauand2023quasi}, the authors established a  connection  between stochastic approximation and ES  by replacing the stochastic perturbations with periodic ones. To minimize  convergence bias, \cite{lauand2023quasi} employs filtering techniques and reduces learning/update gains while keeping the perturbation amplitude constant. The same work also develops state-dependent probing signals to improve transient performance, but still leads to biased convergence.
Regulation of inputs directly to their unknown optimal values is achieved in \cite{grushkovskaya2018class} and \cite{scheinker2014non} using an ES that vanishes at the origin, assuming the optimal value of the cost function is known beforehand. Relaxing this restriction, \cite{abdelgalil2021lie}, \cite{guay2020uncertainty}, \cite{haring2016asymptotic}, and \cite{suttner2019extremum} present control designs that achieve asymptotic convergence of both inputs and outputs to their unknown optimum values. Specifically,  \cite{abdelgalil2021lie} achieves this under certain initial conditions, \cite{haring2016asymptotic} employs time-varying tuning parameters with decaying frequency and amplitude, \cite{guay2020uncertainty} reduces the size of the search region by estimating the uncertainty set around the optimizer and updating the amplitude accordingly, and \cite{suttner2019extremum} adopts an approach with unboundedly growing update rates and frequencies. 

In our prior work \cite{yilmaz2023press}, which is an extended journal version of the conference paper  \cite{yilmaz2023exponential}, we introduce a novel ES called the exponential unbiased extremum seeker (uES), which achieves a first in the field: exponential and unbiased convergence to the unknown optimum at a user-defined rate. However, perfect tracking of time-varying optima remains unaddressed in \cite{yilmaz2023press}. Instead, we discuss the robustness of the design in such scenarios, introducing a modified design with user-adjustable oscillations that gradually decrease but remain non-zero, thus sacrificing perfect convergence.

\subsection{ES for tracking time-varying optima}

Several ES designs have been developed for plants with periodic steady-state outputs and constant optimizers. A scheme in \cite{wang2000extremum} minimizes the limit cycle size by detecting its amplitude and adjusting a controller parameter to reach a constant optimum value. Another design in \cite{guay2007flatness} uses the system's flatness property to track an optimal orbit of a nonlinear dynamical system. For plants with periodic outputs of known periodicity, an ES controller with a moving-average filter is presented in \cite{haring2013extremum}. However, the reliance on periodic performance functions in these designs limits their applicability to a broad and rich set of engineering applications, such as those investigated  in \cite{krstic2014extremum,moura2013lyapunov,wepener2023extremum,hu2015extremum} (and the references therein).

The optimization of arbitrary time-varying cost functions using ES is considered in various studies. The groundwork for time-varying optimizers is laid in \cite{krstic2000performance}, where a generalized ES scheme is developed to track optimizers with known dynamics but unknown coefficients, employing the internal model principle. For slowly varying optima, \cite{sahneh2012extremum} introduces a delay-based strategy to extract the gradient signal. This strategy is further extended to handle dynamic systems with input constraints in \cite{ye2013extremum}.
In \cite{scheinker2012extremum}, an ES is developed to provide unknown reference tracking and stabilization for a class of unknown nonlinear systems, based on time-varying nonlinear high-gain feedback. To address unstable nonlinear systems with time-varying extrema, \cite{moshksar2015model} presents a model-based adaptive ES algorithm. The results in \cite{grushkovskaya2017extremum} and \cite{hazeleger2020extremum} establish  local and semi-global practical asymptotic stability of the extremum with a dynamic map. While \cite{hazeleger2020extremum} seeks a constant optimizer by ES to optimize time-varying steady-state plant, \cite{grushkovskaya2017extremum} extends the Lie bracket approximation method to prove convergence toward a neighborhood of a time-varying optimizer. The robustness of Lie bracket-based ES schemes with respect to time-varying parameters is investigated in \cite{labar2022iss} within the framework of input-to-state stability (ISS). Additionally, \cite{poveda2021fixed} studies ISS-like properties of fixed-time extremum seeking with a time-varying cost function. In contrast to prior results, which utilize high-frequency dither signals, \cite{michael2023gradient} develops a cooperative ES scheme for tracking moving sources without the need for dither signals. Although all these studies achieve convergence to a small neighborhood of the time-varying optimum, they do not ensure perfect tracking. Our recent conference paper \cite{yilmaz2024perfect} presents an ES capable of perfectly tracking an unknown time-varying extremum. However, that design is limited to strongly convex maps and bounded time-varying extrema, unlike the divergent cases considered in this paper.

\subsection{ES with finite, fixed and prescribed-time convergence}
Recent advancements in ES field have introduced concepts such as finite-time ES \cite{guay2021finite} and fixed-time ES \cite{poveda2021fixed}, \cite{poveda2021nonsmooth}, to accelerate convergence to the optimum's vicinity. The key distinction between these concepts is that in finite-time stability, convergence time depends on initial conditions, while in fixed-time stability, it is determined by the system's parameters regardless of initial conditions. Prescribed-time stability, introduced by \cite{song2017time}, strengthens these notions further. It allows users to pre-define a desired terminal time regardless of initial conditions and system’s
parameters. By incorporating this concept into ES, an algorithm called prescribed-time extremum seeking (PT-ES) has been introduced in \cite{todorovski2023practical} and \cite{yilmaz2024prescribed}. However, these designs achieve convergence to a neighborhood of the extremum, not the exact extremum, using unbounded controllers. 
Our companion paper \cite{yilmaz2023press} addresses this issue with an prescribed-time unbiased extremum seeker (PT-uES), which keeps control inputs bounded but is limited to strongly convex maps and a fixed optimum.

\subsection{Methodology: State scaling and time scaling}
Our ES designs are built upon  two main methodological advances: state scaling and time scaling. Based on state scaling, the system states are multiplied by an unboundedly growing function, and the original system is converted to a new ``transformed system''. The stability of the transformed system is equivalent to that of the original system, whether asymptotic, exponential, or within a prescribed time, depending on the growth rate of the scaling function. A crucial aspect of this accelerated control technique is that \emph{the states converge faster than the gains diverge,} which guarantees  boundedness of the input signal.

The state scaling technique has its roots in the mid-60s, with work by \cite{sandberg1965some} and \cite{zames1965nonlinear} on analysis of nonlinear feedback systems, which later became known as the ``exponential weighting technique''. For systems lacking persistency of excitation (PE), exponentially growing gains is used in \cite{nihtila1987exponential} to compensate for the PE loss and achieve unbiased parameter identification at a desired exponential rate. A major advancement in stability is introduced in \cite{song2017time} with a scaling function that grows and blows up at a user-defined time. This has spurred significant interest in prescribed-time control \cite{song2023prescribed}.

Time scaling provides an alternative approach to achieve faster convergence. By conceptually ``dilating'' the fixed time interval to infinity, controller design and stability analysis are performed in a virtual infinite time domain. Crucially, stability proven in this dilated time domain translates back to guaranteed stability over the original finite time interval after ``contracting'' time back.
In this paper, we use this method to dilate time from the domain (not necessarily finite) to a new infinite one, allowing us to handle the arbitrarily varying optimum by slowing its speed in the dilated time domain. This enables tracking of the optimum asymptotically, exponentially, or within a prescribed time.

Throughout this paper, we refer to these techniques as state transformation and time transformation.

\subsection{Contributions}

In this paper, we introduce three uES designs that achieve unbiased convergence to the arbitrarily time-varying optimum: asymptotic uES, exponential uES, and PT-uES, which achieve convergence asymptotically, exponentially, and in prescribed time, respectively. Our key contributions are listed below:
\begin{itemize}
    \item We extend the capabilities of our exponential and PT uES designs from prior work \cite{yilmaz2023press} by $(i)$ broadening the strong convexity assumption to encompass a wider class of convex cost functions, $(ii)$ relaxing the fixed optima assumption to include time-varying optima, and $(iii)$ accommodating $\mathcal{C}^2$ functions (continuously differentiable up to the second order) instead of $\mathcal{C}^4$.
    \item We establish feasibility conditions for selecting time-varying design parameters and their decay/growth rates in relation to the convexity of the map and the decay/growth rate of the optima. These conditions imply that the use of constant-frequency probing restricts the range of achievable asymptotic unbiasing rates, depending on the ``flatness'' of the map at the optimum. In contrast, employing chirpy probing enables users to define the unbiasing rate arbitrarily--whether asymptotic, exponential, or in prescribed time--provided that the frequency and adaptation rates grow sufficiently fast.
    \item Compared to existing designs \cite{abdelgalil2021lie, guay2020uncertainty, haring2016asymptotic, suttner2019extremum}, which achieve asymptotic unbiased convergence, our asymptotic uES design offers convergence at a user-defined asymptotic rate. To the best of our knowledge, perfect tracking of time-varying optima at a desired rate has not been achieved in prior works in the field.
\end{itemize}
Other contributions are as follows:
\begin{itemize}
    \item We develop an exponential uES tailored to a unicycle, which requires singular perturbation analysis in addition to averaging.
    \item We  experimentally test the developed design on a light-seeking task by tuning the angular velocity of a unicycle robot based on real-time light sensor data.
\end{itemize}

This paper is an extended journal version of our conference papers \cite{yilmaz2024unbiased, yilmaz2024perfect} and differs in four key ways: the results are generalized to a broader class of functions, the exponential and PT-uES designs are developed for this class, the optima are allowed to diverge at arbitrary speeds (even in finite time), and experimental results are presented.

\subsection{Organization}
The structure of this paper is as follows. Section \ref{prelim} introduces fundamental stability concepts, which serve as references throughout the paper. Section \ref{problemstate} outlines the problem formulation. Sections \ref{subsec:asymuES}, \ref{expuES}, and \ref{subsec:prescuES} introduce designs for asymptotic, exponential and prescribed-time uES, along with formal stability analysis. Section \ref{sec:experiment} considers the problem of source seeking and presents the experimental results obtained. Finally, Section \ref{sec:conclusion} concludes the paper.

\subsection{Notation}
We denote the 
Euclidean norm by $|\cdot|$. The $\delta$-neighborhood of a set $\mathcal{S} \subset \mathbb{R}^n$ is denoted by $U_{\delta}^{\mathcal{S}}=\{x \in \mathbb{R}^n : \inf_{e \in \mathcal{S}} |x-e|<\delta\}$. The Lie
bracket of two vector fields $f, g : \mathbb{R}^n \times \mathbb{R} \rightarrow \mathbb{R}^n$ with $f (\cdot, t), g( \cdot, t)$ being continuously differentiable is defined by $[f, g](x, t) := \frac{\partial g(x,t)}{\partial x}f(x,t)-\frac{\partial f(x,t)}{\partial x}g(x,t)$. The notation $e_i$ corresponds to the $i$th unit vector in $\mathbb{R}^n$. $\mathbb{R}^+$ denotes the set of non-negative real numbers. We denote the gradient of a function $h: \mathbb{R}^n \to \mathbb{R}$ by $\nabla h(x)$.

\section{Preliminaries} \label{prelim}
\subsection{Lie bracket averaging}
Consider a control-affine system
\begin{align}
    \dot{x}={}f_0(x,t)+\sum_{i=1}^m f_i(x,t) \sqrt{\omega} u_i(\omega t), \label{conaff}
\end{align}
where $x(t_0)=x_0 \in \mathbb{R}^n$, $t_0 \in \mathbb{R}^+$, $\omega>0$, $f_j \in \mathcal{C}^2: \mathbb{R}^n \times \mathbb{R}^+ \to \mathbb{R}^n$ for $j=0,\dots,m$, and $u_i \in \mathcal{C}^1: \mathbb{R}^+  \to \mathbb{R}$ for $i=1,\dots,m$. The functions $u_i$ are $T$-periodic with some $T>0$, and $\int_0^T u_i(\sigma)d\sigma=0$ for $i=1,\dots,m$. We compute the Lie bracket system corresponding to \eqref{conaff} as follows
\begin{align}
    \dot{\bar{x}}={}f_0(\bar{x},t)+\frac{1}{T} \sum_{{i=1}\atop{j=i+1}}^{m} [f_i, f_j](\bar{x},t) \int_0^T \int_0^{\vartheta}  u_{j}(\vartheta) u_{i}(\sigma) d\sigma d\vartheta, \label{conafflie}
\end{align}
where $\bar{x}(t_0)=x(t_0)$. We make the following assumption regarding the boundedness of vector fields and their partial derivatives:
\begin{assumption} \label{vectbound}
For every compact set $\mathcal{X} \subset \mathbb{R}^n$, $|f_i(x,t)|$, $\left| \frac{\partial f_i(x,t)}{\partial t} \right|$, $\left| \frac{\partial f_i(x,t)}{\partial x} \right|$, $\left| \frac{\partial^2 f_j(x,t)}{\partial t \partial x} \right|$, $\left| \frac{\partial [f_j, f_k](x,t)}{\partial t} \right|$, $\left| \frac{\partial [f_j, f_k](x,t)}{\partial x} \right|$ are bounded for all $x \in \mathcal{X}$, $t \geq t_0$, $i=0,\dots,m$, $j=1,\dots,m$, $k=j,\dots,m$.    
\end{assumption}

To understand the relationship between the stability of \eqref{conaff} and \eqref{conafflie}, we revisit a theorem presented in \cite{durr2013lie}. 
Refer to Definition \ref{def:localpractical} in Appendix \ref{sec:app2} for the formal definition of practical stability.

\begin{theorem}[\cite{durr2013lie}]   \label{LieBracketAvThe}
Consider the system \eqref{conaff} under Assumption \ref{vectbound}.
If a compact set $\mathcal{S} \subset \mathbb{R}^n$ is locally (globally) uniformly asymptotically stable for \eqref{conafflie}, then $\mathcal{S}$ is locally (semi-globally) practically uniformly asymptotically stable for \eqref{conaff}. 
\end{theorem}

\subsection{Singularly perturbed Lie bracket averaging} \label{tit:sing_pert}
Consider the system of the form
\begin{align}
    \dot{x}={}&\varepsilon f_0(x, z, \varepsilon t)+\varepsilon \sqrt{\omega} \sum_{i=1}^m f_i(x,z,\varepsilon t)  u_i(\omega \varepsilon t), \label{singpert1}  \\
    \dot{z}={}&g(x,z), \label{singpert2}
\end{align}
with $x(t_0)=x_0 \in \mathbb{R}^n$, $z(t_0)=z_0 \in \mathbb{R}^m$, $t_0 \in \mathbb{R}^+$, $\varepsilon, \omega>0$, $f_j \in \mathcal{C}^2: \mathbb{R}^n \times \mathbb{R}^m \times \mathbb{R}^+ \to \mathbb{R}^n$ for $j=0,\dots,m$, and $u_i \in \mathcal{C}^1: \mathbb{R}^+  \to \mathbb{R}$ for $i=1,\dots,m$. The functions $u_i$ are $T$-periodic with some $T>0$, and $\int_0^T u_i(\sigma)d\sigma=0$ for $i=1,\dots,m$.
To analyze stability, we first find the 
quasi-steady-state, denoted by $l(x)$, which satisfies $g(x, l(x)) = 0$. Performing the change of variables $z_{b} = z-l(x)$, we get the boundary layer model as
\begin{align}
    \dot{z}_b={}g(x,z_b+l(x)). \label{bound_sing}
\end{align}
Then, by substituting the quasi-steady state into \eqref{singpert1}, we obtain a reduced model as
\begin{align}
    \dot{x}_r={}\varepsilon f_0(x_r, l(x_r), \varepsilon t)+\varepsilon \sqrt{\omega} \sum_{i=1}^m f_i(x_r,l(x_r),\varepsilon t)  u_i(\omega \varepsilon t),
\end{align}
in which the corresponding Lie bracket average, for $\varepsilon=1$, is given by
\begin{align}
    \dot{\bar{x}}_r={}& f_0(\bar{x}_r, l(\bar{x}_r),  t)+ \frac{1}{T} \sum_{{i=1}\atop{j=i+1}}^{m} [f_i, f_j](\bar{x}_r, l(\bar{x}_r), t) \nonumber \\
    &\times \int_0^T \int_0^{\vartheta}  u_{j}(\vartheta) u_{i}(\sigma) d\sigma d\vartheta. \label{avredusing}
\end{align}
We make the following assumption, similar to Assumption \ref{vectbound}:
\begin{assumption}  \label{ass:singular}
For every compact sets $\mathcal{X} \subset \mathbb{R}^n$ and $\mathcal{Z} \subset \mathbb{R}^m$, $|f_i(x,z,t)|$, $\left| \frac{\partial f_i(x,z,t)}{\partial t} \right|$, $\left| \frac{\partial f_i(x,z,t)}{\partial x} \right|$, $\left| \frac{\partial f_i(x,z,t)}{\partial z} \right|$, $\left| \frac{\partial^2 f_j(x,z,t)}{\partial t \partial x}\right|$, $\left| \frac{\partial^2 f_j(x,z,t)}{\partial t \partial z}\right|$, $\left| \frac{\partial [f_j, f_k](x,z,t)}{\partial t} \right|$, $\left| \frac{\partial [f_j, f_k](x,z,t)}{\partial x} \right|$, $\left| \frac{\partial [f_j, f_k](x,z,t)}{\partial z} \right|$ are bounded for all $x \in \mathcal{X}$, $z \in \mathcal{Z}$, $t \geq t_0$, $i=0,\dots,m$, $j=1,\dots,m$, $k=j,\dots,m$.    
\end{assumption}
We present the following theorem from \cite{durr2017extremum} and provide the corresponding stability definition in Definition \ref{def:singular} in Appendix \ref{sec:app2}.
\begin{theorem}[\cite{durr2017extremum}] \label{thm:singperturb}
Consider the system \eqref{singpert1} and \eqref{singpert2} under Assumption \ref{ass:singular}. Suppose that a compact set  $\mathcal{S} \subset \mathbb{R}^n$ is globally uniformly asymptotically stable for the average of the reduced system \eqref{avredusing} and the origin of the boundary layer model \eqref{bound_sing} is globally exponentially stable. Then, the set $\mathcal{S}$ is singularly semi-globally practically uniformly asymptotically stable for \eqref{singpert1}. 
\end{theorem}

\section{Problem Statement} \label{problemstate}

We consider the following optimization problem
\begin{align}
    \min_{\theta \in \, \mathbb{R}^n} J(\theta), \label{Jopt}
\end{align}
where $\theta \in \mathbb{R}^n$ is the input, $J \in \mathcal{C}^2: \mathbb{R}^n \to \mathbb{R}$
is an unknown cost function. We make the following assumptions regarding the unknown static map $J(\cdot)$. 
\begin{assumption} \label{Ass0}
    The cost function $J(\cdot)$ has a unique minimum at $\theta^*(t) \in \mathbb{R}^n$, i.e., $J(\theta^*(t)) < J(\theta)$ for $\forall \theta \neq \theta^*(t)$ with $\nabla J(\theta^*(t))=0$.
\end{assumption}
\begin{assumption} \label{Ass1}
There exist constants $\rho_1$, $\rho_2$, $\rho_3$, $\rho_4>0$ and $\kappa \in \mathbb{N}$ such that
\begin{align}
     &\rho_1 |\theta-\theta^*(t)|^{2 \kappa } \leq {}  J(\theta)-J(\theta^*(t)) \leq \rho_2|\theta-\theta^*(t)|^{2 \kappa}, \label{Jboundcon} \\
    &(\theta-\theta^*(t))^T \nabla J(\theta)
    \leq{} \rho_3 |\theta-\theta^*(t)|^{2\kappa}, \\
    &\left|\frac{\partial^2 J(\theta)}{\partial \theta^2}\right| \leq{} \rho_4  |\theta-\theta^*(t)|^{2 \kappa -2} \label{Jboundconen}
\end{align}
\end{assumption}
for all $\theta, \theta^*(t) \in \mathbb{R}^n$.
\begin{assumption} \label{asympextbound}
The derivatives of the time-varying optimizer ${\theta}^*(t)$ and the time-varying optimum $J(\theta^*(t))$ up to the second order obey the following bounds:
\begin{align}
|\dot{\theta}^*(t)|+|\ddot{\theta}^*(t)| \leq {} & M_{\theta} \phi^{c}(t), \quad t \in [t_0, \infty), \label{ass3bound} \\
\left| \dot{J}(\theta^*(t))\right|+\left| \ddot{J}(\theta^*(t))\right| \leq {}&  M_J\phi^{d}(t),   \quad t \in [t_0, \infty),\label{ass3bound2} 
\end{align}
with unknown constants $M_{\theta}, M_{J} \geq  0$, $c, d \in \mathbb{R}$, and a continuous strictly increasing function $\phi(t): \mathbb{R}^+ \to \mathbb{R}$. 
\end{assumption}

Assumption \ref{Ass0} guarantees the existence of a minimum of the function $J(\theta)$ at $\theta=\theta^*(t)$. 
Assumption \ref{Ass1} establishes  bounds on the cost function $J(\theta)$, its gradient, and its Hessian.
These bounds exhibit specific growth rates governed by power functions, similar to those exploited in \cite{grushkovskaya2018class} and \cite{todorovski2023practical}.
For $\kappa=1$, this simplifies to a more conservative assumption of strong convexity, which is considered in \cite{yilmaz2023press}. For $\kappa>1$, the assumption requires a specific class of convexity for $J(\theta)$. Assumption \ref{asympextbound} imposes bounds on the growth or decay of the time-varying optimizer and optimum, based on the parameters $c$ and $d$. Note that this assumption encompasses scenarios where the optimum diverges asymptotically, exponentially, or even in finite time.

Using the first order condition of convexity \cite[pp. 69]{boyd2004convex}, we obtain from \eqref{Jboundcon} that
\begin{align}
    (\theta-\theta^*(t))^T \nabla J(\theta)\geq{}& \rho_1 |\theta-\theta^*(t)|^{2\kappa}. \label{Jconvexityfirst}
\end{align}

We measure the unknown function $J(\theta)$ in real time as 
\begin{align}
    y(t)={}&J(\theta(t)), \qquad t \in [t_0,\infty),
\end{align}
in which $y \in \mathbb{R}$ is the output.
Our aim is to design ES algorithms 
using output feedback $y(t)$ in order to achieve the unbiased convergence of $\theta$ to $\theta^{*}(t)$ while simultaneously minimizing the steady state value of $y$, without requiring prior knowledge of either the optimum input $\theta^*(t)$ or the function $J(\cdot)$.

\begin{figure}[t]
    \centering
    \includegraphics[width=0.9\columnwidth]{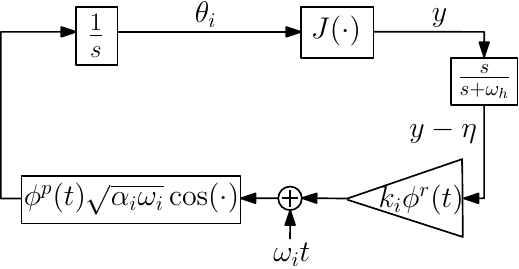} 
    \caption{Unbiased ES scheme with constant-frequency probing for the $i$th element $\theta_i$ of $\theta$, which guarantees $\theta \to \theta^*(t)$ at the rate that $1/\phi(t)$ converges to zero. Refer to Table \ref{table_consfreq} for the function and parameters. }
    \label{uES_consfreq}
\end{figure}

\begin{figure}[t]
    \centering
    \includegraphics[width=\columnwidth]{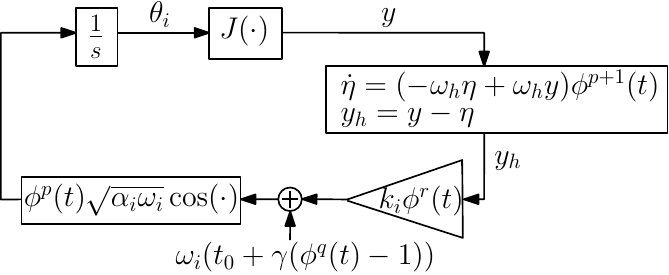} 
    \caption{Unbiased ES scheme with chirpy probing for the $i$th element $\theta_i$ of $\theta$, which guarantees $\theta \to \theta^*(t)$ at the rate that $1/\phi(t)$ converges to zero. Refer to Table \ref{table_chirpy} for the function and parameters.}
    \label{uES_chirpy}
\end{figure}

In the subsequent sections, we present three types of unbiased ES algorithms: asymptotic uES, exponential uES, and prescribed-time uES. These algorithms aim to achieve the aforementioned objective asymptotically, exponentially, and in user-prescribed time. 
For visualization, Figures \ref{uES_consfreq} and \ref{uES_chirpy} depict the uES designs schematically. These figures illustrate constant-frequency probing and chirpy probing approaches, which can handle both convergent and divergent optimums.
Further details on amplitude and gain functions can be found in Tables \ref{table_consfreq} and \ref{table_chirpy}. It is important to note that our uES design modifies the alternative ES design \cite{scheinker2014extremum}, which can only ensure practical stability around the fixed optimum, by incorporating time-varying design parameters.

\begin{remark}
In the unbiased ES designs depicted in Fig. \ref{uES_consfreq}, the gain $k_i \phi^r(t)$ experiences either asymptotic or exponential growth based on the chosen ES type, while the update rate $\phi^p(t) \sqrt{\alpha_i \omega_i}$ experiences asymptotic or exponential decay corresponding to the specific ES type. The crucial aspect of our designs is that the high-pass filtered state $y-\eta$ decays to zero at least as fast as the gain $k_i \phi^r(t)$ grows, ensuring the boundedness of the resulting signal. For a more detailed analysis, refer to Theorem \ref{theoremasymp} and \ref{theoremexp}.
\end{remark}

\begin{table}[t]
\scalebox{1}{
\begin{tabular}{|l|c|c|l|}
\hline
                  & $\phi(t)$   &  $p$   & Conditions        \\ \hline
 & \multirow{5}{*}{$(1+\beta \ (t-t_0))^{\frac{1}{v}}$} & \multirow{5}{*}{$-1$}  & $\kappa \geq 1$  \\ 
 & &   & $v>2\kappa-r \geq 0$ \\ 
\textbf{Asymptotic}& & & $k_i, \alpha_i, \beta>0$ \\
 \textbf{uES} & &  & $c<-1-2\kappa+r$ \\ 
 & &   & $d<-2\kappa$ \\ \hline
 & \multirow{6}{*}{$e^{\lambda(t-t_0)}$} & \multirow{6}{*}{$-1$}  & $\kappa=1$ \\
  &  &  & $r=2$  \\ 
\textbf{Exponential} & & & $k_i \alpha_i>2\lambda/\rho_1$\\
\textbf{uES } &  &  & $0<\lambda<\omega_h/2$  \\ 
  & & & $c<-1$\\
 & &   & $d<-2$ \\\hline
\end{tabular}} 
\caption{Time-varying functions used in Figure \ref{uES_consfreq}.}
\label{table_consfreq}
\end{table}

\begin{remark} \label{remark:rv}
Another point that requires discussion is how parameters $r$ and $v$ influence the performance of asymptotic uES with constant-frequency probing in Fig. \ref{uES_consfreq}.
The constraint on $r$ is $r\leq 2\kappa$, as provided in Table \ref{table_consfreq}. When $\kappa$ is known, one can choose $r=2\kappa$ and satisfy the condition $v>2\kappa-r \geq 0$ with any arbitrarily small $v$. This makes the rate of growth of $\phi(t)$ faster, which in turn accelerates convergence. For a $\kappa$-blind user, the ``safe'' choice of $r$ and $v$ would be 2 and a sufficiently large value, respectively, compromising the rate of convergence.
\end{remark}

\begin{remark}
The uES designs with chirpy probing in Table~\ref{table_chirpy} relax the conditions on $c$ and $d$ in \eqref{ass3bound} and \eqref{ass3bound2}, which characterize the rate of growth of optima. A divergent optimum with any growth rate can be tracked perfectly with a sufficiently large $p$, which determines the growth rate of the adaptation and frequency. In addition, compared to the exponential uES in Table~\ref{table_consfreq}, the exponential uES with chirp probing in Table~\ref{table_chirpy} achieves convergence for any $\kappa \geq 1$.
\end{remark}

\begin{remark} \label{remark:saturation}
In practical scenarios, the implementation of the developed algorithms might be limited by various factors. To mitigate this, the time-varying signal $\phi(t)$ along with the instantaneous frequency $\omega_i \gamma d\phi^q/dt$ in the case of chirped probing, can be constrained to moderately large values that are sufficient for close tracking of the extremum.
\end{remark}

\begin{table}[t]
\scalebox{0.98}{
\begin{tabular}{|l|c|c|l|}
\hline
                  & $\phi(t)$   &  ${p}$   & Conditions        \\ \hline
& \multirow{5}{*}{$(1+\beta \ (t-t_0))^{\frac{1}{v}}$} & \multirow{5}{*}{$q-v-1$}  & $\kappa \geq 1$  \\ 
 & &   & ${q}>2\kappa-r\geq 0$ \\ 
\textbf{Asymptotic} & &   & $\gamma=v/(\beta q)$ \\ 
\textbf{uES} & &  & $c<p-2\kappa+r$ \\ 
 & &   & $d<p-2\kappa+1$ \\ \hline
 & \multirow{5}{*}{$e^{\lambda(t-t_0)}$} & \multirow{5}{*}{$q-1$}  & $\kappa \geq 1$ \\
 &  &  & ${q}>2\kappa-r\geq 0$  \\ 
 \textbf{Exponential}&  &  & $\gamma=1/(\lambda q)$  \\ 
 \textbf{uES }&  &  & $c<p-2\kappa+r$  \\ 
 & &   & $d<p-2\kappa+1$ \\\hline
 & \multirow{5}{*}{$\left(\frac{T}{T+t_0-t}\right)^{\frac{1}{\varrho}}$} & \multirow{5}{*}{$q+\varrho-1$}  & $\kappa \geq 1$ \\
 &  &  & ${q}>2\kappa-r\geq 0$  \\ 
 \textbf{Prescribed-time}&  &  & $\gamma=\varrho T/q$  \\ 
 \textbf{uES} &  &  & $c<p-2\kappa+r$  \\ 
 & &   & $d<p-2\kappa+1$ \\\hline
\end{tabular}} \\ \\
\caption{Time-varying functions used in Figure \ref{uES_chirpy} with any $\beta, \lambda, v, \varrho, k_i, \alpha_i >0$ for $i=1,\dots,n$.}
\label{table_chirpy}
\end{table}


\section{Asymptotic uES} \label{subsec:asymuES}

\subsection{Asymptotic uES with constant-frequency probing}
We introduce an ES with asymptotic and unbiased convergence in the following theorem.
\begin{theorem} \label{theoremasymp}
Consider the following asymptotic uES design
\begin{equation} \label{ESasymp}
\begin{cases}
\begin{aligned}
\dot{\theta}={}&\xi^{-1}(t) \sum_{i=1}^n \sqrt{\alpha_i \omega_i}  e_i \cos\Big(\omega_i t+k_i \xi^r(t)(J(\theta)-\eta)\Big),   \\
\dot{\eta}={}&-\omega_h \eta+\omega_h J(\theta),
\end{aligned}
\end{cases}
\end{equation}
with the asymptotically growing function
\begin{align}
    \xi(t)={}(1+\beta(t-t_0))^{\frac{1}{v}}, \qquad t \in [t_0, \infty). \label{xidotdyn}
\end{align}
Let $\omega_i = \omega \hat{\omega}_i$ such that $\hat{\omega}_i \neq \hat{\omega}_j$ $\forall i \neq j$, $ t_0 \geq 0$, $\alpha_i, k_i, \beta, \omega_h>0$ $\forall i=1,\dots,n$ and $v > 2\kappa-r \geq 0$.  Let Assumptions \ref{Ass0}, \ref{Ass1} hold, and Assumption \ref{asympextbound} holds with $\phi(t)=\xi(t)$, $c<-1-2\kappa+r$, $d<-2\kappa$. Then, there exists $\omega^*>0$ such that for all $\omega > \omega^*$, 
the following holds
\begin{itemize}
    \item $\theta(t) \to \theta^*(t)$ semi-globally with respect to $\omega$ and asymptotically at the rate of $1/\xi(t)$, and there exist a class $\mathcal{KL}$ function $\mathcal{B}$, a class $\mathcal{K}$ function $\mathcal{Y}$, and a nonnegative constant $D(\theta(t_0), \theta^*(t_0), \eta(t_0))$ such that 
    \begin{align}
    |\theta(t)-\theta^*(t)|&\leq{}\xi^{-1}(t) \Big(D+\mathcal{B}\left(|\theta(t_0)-\theta^*(t_0)|, t-t_0\right) \nonumber \\
    &\hspace{0.3cm}+\sup_{t_0 \leq s \leq t}\mathcal{B}\left(\mathcal{Y}\left( \xi^{-1}(s) \right), t-s\right)\Big). \label{KLforasymp}
\end{align}
    \item $\eta(t), y(t) \to J(\theta^*(t))$ semi-globally with respect to $\omega$ and asymptotically  at the rate of $1/\xi^{2\kappa}(t)$.
\end{itemize}
\end{theorem}

\begin{proof} Let us proceed through the proof step by step.

\textbf{Step 1: State transformation.}
Let us consider the following transformations
\begin{align}
    \theta_{f}={}& \xi(t) (\theta-\theta^*(t)), \label{statetransfor} \\
    \eta_{f}={}& \xi^{2\kappa}(t) (\eta-J(\theta^*(t))), \label{statetransfor2}
\end{align}
which transform \eqref{ESasymp}  to
\begin{equation} \label{thetaasymtrans}
\begin{cases}
\begin{aligned}
\dot{\theta}_{f}={}&-\xi(t)\dot{\theta}^*(t)+\frac{\beta}{v } \xi^{-v}(t) {\theta}_{f}+ \sum_{i=1}^n \sqrt{\alpha_i \omega_i}  e_i  \\
&\times \cos\bigg(\omega_i t+k_i \xi^r(t)\left(J_f(\theta_f, t)-\xi^{-2\kappa}(t){\eta}_f\right)\bigg),  \qquad \quad \\
\dot{\eta}_f={}&\left(\frac{2\beta \kappa}{v}\xi^{-v}(t)-\omega_h \right)\eta_f+\omega_h \xi^{2\kappa}(t) J_f(\theta_f, t)  \\
&-\xi^{2\kappa}(t)\dot{J}(\theta^*(t)), 
\end{aligned}
\end{cases}
\end{equation}
with
\begin{align}  
    J_f(\theta_f, t) = {}& J(\theta_f/\xi(t)+\theta^*(t)) -J(\theta^*(t)).
\end{align}
The initial and most critical step in our analysis relies on transformations \eqref{statetransfor}, \eqref{statetransfor2}. These establish a precondition:  once we demonstrate the stability of the transformed system \eqref{thetaasymtrans}, the asymptotic convergence of $\theta$ to $\theta^*(t)$ and $\eta$ to $J(\theta^*(t))$ naturally follows. To analyze the stability of \eqref{thetaasymtrans}, we employ the Lie bracket averaging technique. As a first step, we rewrite the $\theta_f$-system in \eqref{thetaasymtrans} by expanding the cosine term, as shown below
\begin{align}
    \dot{\theta}_{f}=&-\xi(t)\dot{\theta}^*(t)+\frac{\beta}{v} \xi^{-v}(t) {\theta}_{f}+ \sum_{i=1}^n \sqrt{\alpha_i \omega_i} e_i \cos(\omega_i t) \nonumber \\
    &\times \cos\left(k_i \xi^r(t) \left(J_f(\theta_f, t)-\xi^{-2\kappa}(t){\eta}_f\right)\right)-\sum_{i=1}^n \sqrt{\alpha_i \omega_i}\nonumber \\ 
    & \times e_i \sin(\omega_i t)  \sin\left(k_i \xi^r(t) \left(J_f(\theta_f, t)-\xi^{-2\kappa}(t){\eta}_f\right)\right). \label{thetaasymtranse} 
\end{align}
Next, we integrate the $\eta_f$-system from \eqref{thetaasymtrans} with the transformed $\theta_f$-system \eqref{thetaasymtranse} to construct the following system 
\begin{align}
    \begin{bmatrix}
        \dot{\theta}_f \\ \dot{\eta}_f
    \end{bmatrix}={}&b_0(\theta_f, \eta_f,t)+\sum_{i=1}^n b_{c,i}(\theta_f, \eta_f,t) \sqrt{\omega_i} \cos(\omega_i t) \nonumber \\
    &-b_{s,i}(\theta_f, \eta_f,t) \sqrt{\omega_i} \sin(\omega_i t), \label{transfcomp}
\end{align}
where
\begin{align}
    &b_0(\theta_f, \eta_f,t)\nonumber \\
    &=\begin{bmatrix}
        -\xi(t)\dot{\theta}^*(t)+\frac{\beta}{v \xi^{v}(t)}  \theta_f \\
        \left(\frac{2\beta \kappa}{v \xi^{v}(t)} -\omega_h \right)\eta_f +\xi^{2\kappa}(t) \left(\omega_h  J_f(\theta_f, t) -\dot{J}(\theta^*(t))\right)
    \end{bmatrix}, \label{bodefin} \\
    &b_{c,i}(\theta_f, \eta_f,t)\nonumber \\
    &=\begin{bmatrix} \sqrt{\alpha_i} e_i \cos\Big(k_i \xi^r(t) \left(J_f(\theta_f, t)- 
     \xi^{-2\kappa} (t){\eta}_f \right)\Big) \\ 0 \end{bmatrix}, \\
    &b_{s,i}(\theta_f, \eta_f,t) \nonumber \\
    &=\begin{bmatrix} \sqrt{\alpha_i} e_i \sin\Big(k_i \xi^r(t) \left(J_f(\theta_f, t)-\xi^{-2\kappa} (t){\eta}_f\right)\Big) \\ 0 \end{bmatrix}.
\end{align}

\textbf{Step 2: Feasibility analysis of \eqref{transfcomp} for averaging.}
To ensure the applicability of Lie bracket averaging, we need to verify that \eqref{transfcomp} satisfies the boundedness assumption outlined in Theorem \ref{LieBracketAvThe}. Towards this end, let $\mathcal{M} \subset \mathbb{R}^n$ and $\mathcal{Z} \subset \mathbb{R}$ be compact sets, and consider the bounds \eqref{Jboundcon}, \eqref{ass3bound}, and \eqref{ass3bound2}. Then, we derive 
\begin{align}
    \xi^{2\kappa}(t)  |J_f(\theta_f, t)| \leq{}& \rho_2 |\theta_f|^{2\kappa}, \label{Jfxibound} \\
     \xi^{2\kappa}(t) \left( |\dot{J}(\theta^*(t))|+|\ddot{J}(\theta^*(t))| \right) \leq{}&M_J \xi^{d+2\kappa}(t), \\
    \xi(t) \left(|\dot{\theta}^*(t)|+|\ddot{\theta}^*(t)| \right) \leq{}& M_{\theta} \xi^{c+1}(t). \label{thetaxibound}
\end{align}
Note that $c+1<0$ due to $c<-1-2\kappa+r$ and $r \leq 2\kappa$, and also that $d+2\kappa<0$. Based the bounds \eqref{Jfxibound}--\eqref{thetaxibound}, we establish the boundedness of $|b_0(\theta_f,\eta_f,t)|$, $|b_{c,i}(\theta_f,\eta_f,t)|$, and $|b_{s,i}(\theta_f,\eta_f,t)|$ for $(t,\theta_f, \eta_f) \in  [t_0, \infty) \times \mathcal{M} \times \mathcal{Z}$. Noting from \eqref{statetransfor} that
\begin{align}
    \frac{d \theta}{d \theta_f}={}\frac{1}{\xi(t)},
\end{align}
and recalling the bounds in Assumption \ref{Ass1}, we derive the following inequalities
\begin{flalign}
    &\xi^{2\kappa}(t) \left|\frac{\partial J_f(\theta_f, t)}{\partial \theta_f}   \right| ={}\xi^{2\kappa-1}(t) \left|\frac{\partial J(\theta)}{\partial \theta}   \right| \leq{}\rho_3 |\theta_f|^{2\kappa-1}, \label{xidelJf1} \\ 
    &\xi^{2\kappa}(t) \left|\frac{\partial^2 J_f(\theta_f, t)}{\partial \theta_f^2}   \right| ={}\xi^{2\kappa-2}(t) \left|\frac{\partial^2 J(\theta)}{\partial \theta^2}   \right| \leq{}\rho_4 |\theta_f|^{2\kappa-2},  \\
    &\xi^{2\kappa}(t) \left|\frac{\partial J_f(\theta_f, t)}{\partial t}\right| \leq {} \xi^{2\kappa}(t)\Bigg( |\dot{J}(\theta^*(t))|+\bigg|\frac{\partial J(\theta)}{\partial \theta} \bigg| \nonumber \\
    &\hspace{3.15cm}\times \left(|\theta_f |\frac{\dot{\xi}(t)}{\xi^2(t)}+|\dot{\theta}^*(t)|\right)  \Bigg), \nonumber \\
    &\hspace{2.75cm} \leq{} M_J \xi^{d+2\kappa}(t)+\rho_3 |\theta_f|^{2\kappa-1} \nonumber \\
    &\hspace{3cm}\times (\beta \xi^{-v}(t) |\bar{\theta}_f|/v+M_{\theta}\xi^{1+c}(t)), \label{xidelJf2} \\
    &\xi^{2\kappa}(t) \left|\frac{\partial^2 J_f(\theta_f, t)}{\partial \theta_f \partial t}   \right| =
    {}\xi^{2\kappa}(t) \left| \frac{\partial}{\partial t}\left( \frac{\partial J(\theta)}{\partial \theta} \frac{1}{\xi(t)}  \right) \right|, \nonumber \\
    &\hspace{2.4cm}\leq{}\xi^{2\kappa-1}(t)  \left| \frac{\partial^2 J(\theta)}{\partial \theta^2} \right| \Bigg( |\theta_f| \frac{\dot{\xi}(t)}{\xi^2(t)}+|\dot{\theta}^*(t)|\Bigg) \nonumber \\
    &\hspace{2.8cm}+\xi^{2\kappa}(t)\left| \frac{\partial J(\theta)}{\partial \theta} \right| \frac{\dot{\xi}(t)}{\xi^2(t)} \nonumber \\
     &\hspace{2.4cm}\leq{} (\beta \rho_3/v) \xi^{-v}(t) |\theta_f|^{2\kappa-1}+\rho_4  |\theta_f|^{2\kappa-2}\nonumber \\
     &\hspace{2.85cm}\times (\beta \xi^{-v}(t)|\theta_f|/v +M_{\theta}\xi^{c+1}(t)),  \label{xidelJf3}
\end{flalign}
which are bounded for $(t,\theta_f, \eta_f) \in  [t_0, \infty) \times \mathcal{M} \times \mathcal{Z}$. 
Considering the bounds \eqref{Jfxibound}--\eqref{thetaxibound}, \eqref{xidelJf1}--\eqref{xidelJf3} and conditions in Table \ref{table_consfreq}, we establish the boundedness of $\left| \frac{\partial b_0(\theta_f,\eta_f,t)}{\partial \theta_f} \right|$, $\left| \frac{\partial b_0(\theta_f,\eta_f,t)}{\partial \eta_f} \right|$, $\left| \frac{\partial b_0(\theta_f,\eta_f,t)}{\partial t} \right|$, $\left| \frac{\partial b_{c,i}(\theta_f,\eta_f,t)}{\partial \theta_f} \right|$, $\left| \frac{\partial b_{s,i}(\theta_f,\eta_f,t)}{\partial \theta_f} \right|$, $\left| \frac{\partial b_{c,i}(\theta_f,\eta_f,t)}{\partial \eta_f} \right|$, $\left| \frac{\partial b_{s,i}(\theta_f,\eta_f,t)}{\partial \eta_f} \right|$, $\left| \frac{\partial b_{c,i}(\theta_f,\eta_f,t)}{\partial t} \right|$, $\left| \frac{\partial b_{s,i}(\theta_f,\eta_f,t)}{\partial t} \right|$, $\left| \frac{\partial^2 b_{s,i}(\theta_f,\eta_f,t)}{\partial \theta_f \partial t} \right|$, $\left| \frac{\partial^2 b_{c,i}(\theta_f,\eta_f,t)}{\partial \theta_f \partial t} \right|$, $\left| \frac{\partial^2 b_{s,i}(\theta_f,\eta_f,t)}{\partial \eta_f \partial t} \right|$, and $\left| \frac{\partial^2 b_{c,i}(\theta_f,\eta_f,t)}{\partial \eta_f \partial t} \right|$ within the domain $(t,\theta_f, \eta_f) \in  [t_0, \infty) \times \mathcal{M} \times \mathcal{Z}$. Next,  we compute the following Lie bracket
\begin{align}
    &\begin{bmatrix}
        b_{c,i}(\theta_f,\eta_f, t) & b_{s,i}(\theta_f,\eta_f, t)
    \end{bmatrix} \nonumber \\
    &\hspace{0.55cm} ={}k_i \sqrt{\alpha_i} \xi^r(t) \begin{bmatrix}   e_i  \frac{\partial J_{f}(\theta_f, t)}{\partial {\theta}_{f}}  \cos(\rho) &  -e_i \xi^{-2\kappa}(t) \cos(\rho) \\ 0 & 0 \end{bmatrix} \nonumber \\
    &\hspace{0.9cm}\times b_{c,i}(\theta_f, \eta_f,t) \nonumber \\
    &\hspace{0.9cm} +k_i \sqrt{\alpha_i} \xi^r(t) \begin{bmatrix}  e_i \frac{\partial J_{f}(\theta_f, t)}{\partial {\theta}_{f}} \sin(\rho) & -e_i \xi^{-2\kappa}(t)  \sin(\rho) \\ 0 & 0 \end{bmatrix} \nonumber \\
    &\hspace{0.9cm}\times b_{s,i}(\theta_f,\eta_f, t), \nonumber \\
    &\hspace{0.55cm}={}k_i \alpha_i \xi^r(t) \begin{bmatrix}e_i \frac{\partial J_{f}(\theta_f, t)}{\partial {\theta}_{f,i}} \\ 0 \end{bmatrix}, \label{liebcbs}
\end{align}
where
\begin{align}
    \rho={}k_i \xi^r(t) \left(J_f(\theta_f, t)-\xi^{-2\kappa}(t){\eta}_f\right).
\end{align}
The boundedness of $\left|\frac{\partial [b_{c,i}(\theta_f,\eta_f, t), \, b_{s,i}(\theta_f,\eta_f, t)]}{\partial \theta_{f}}\right|$ and $\left|\frac{\partial [b_{c,i}(\theta_f, \eta_f, t), \, b_{s,i}(\theta, \eta_f, t)]}{\partial t}\right|$ for $r \leq 2\kappa $ in  $(t,\theta_f, \eta_f) \in  [t_0, \infty) \times \mathcal{M} \times \mathcal{Z}$ is established by recalling the bounds \eqref{xidelJf1}--\eqref{xidelJf3}. Additionally, it should be noted that $\left|\frac{\partial [b_{c,i}(\theta_f,\eta_f, t), \, b_{s,i}(\theta_f,\eta_f, t)]}{\partial \eta_{f}}\right|=0$. As a result, we fulfill the boundedness requirement for Lie bracket averaging. 

\textbf{Step 3: Lie bracket averaging.} We derive the Lie bracket system for \eqref{transfcomp} as follows
\begin{align}
    \begin{bmatrix}
        \dot{\bar{\theta}}_f \\ \dot{\bar{\eta}}_f
    \end{bmatrix}={}&b_0(\bar{\theta}_f, \bar{\eta}_f,t)- \frac{1}{2} \sum_{i=1}^n \begin{bmatrix}
        b_{c,i}(\bar{\theta}_f,\bar{\eta}_f, t) & b_{s,i}(\bar{\theta}_f,\bar{\eta}_f, t)
    \end{bmatrix}. \label{thetabarLie}
\end{align}
Considering \eqref{bodefin} and \eqref{liebcbs}, we can express \eqref{thetabarLie} as
\begin{equation} \label{thetaasymlieball}
\begin{cases}
\begin{aligned}
    \dot{\bar{\theta}}_f={}&-\xi(t)\dot{\theta}^*(t)+\frac{\beta}{v} \xi^{-v}(t) \bar{\theta}_f-\sum_{i=1}^n \frac{k_i \alpha_i}{2}  \xi^{r}(t) e_i \frac{\partial J_f(\bar{\theta}_f, t)}{\partial \bar{\theta}_{f,i}},  \\
    \dot{\bar{\eta}}_f={}&\left(\frac{2\beta \kappa}{v}\xi^{-v}(t)-\omega_h \right)\bar{\eta}_f+\omega_h \xi^{2\kappa}(t) J_f(\bar{\theta}_f, t) \\
    &-\xi^{2\kappa}(t)\dot{J}(\theta^*(t)). 
\end{aligned}
\end{cases}
\end{equation}

\textbf{Step 4: Stability of average system.}
Let us consider the following Lyapunov function
\begin{align}
    V(\bar{\theta}_f)=\frac{1}{2}|\bar{\theta}_f|^2. \label{lyapasymp}
\end{align}
Considering Assumption \ref{Ass1} and \eqref{Jconvexityfirst}, we compute the time derivative of \eqref{lyapasymp} along with \eqref{thetaasymlieball} as
\begin{align}
    \dot{V} \leq{} & \xi(t)|\bar{\theta}_f| |\dot{\theta}^*(t)|+\frac{\beta}{v} \xi^{-v}(t) |\bar{\theta}_f|^2-\sum_{i=1}^n \frac{k_i \alpha_i}{2}  \xi^{r}(t)  \bar{\theta}_{f,i} \nonumber \\
    &\times \frac{\partial J_f(\bar{\theta}_f, t)}{\partial \bar{\theta}_{f,i}}, \nonumber \\
    \leq{} & M_{\theta}\xi^{1+c}(t)|\bar{\theta}_f| +\frac{\beta}{v}\xi^{-v}(t)|\bar{\theta}_f|^2-\frac{ (k \alpha)_{\text{min}}\rho_1}{2}\xi^{-2\kappa+r}(t) \nonumber \\
    &\times|\bar{\theta}_f|^{2\kappa}, \label{lyapaverasymp}
\end{align}
where $(k \alpha)_{\text{min}}=\min\{k_i\alpha_i\}$ for $i=1, \dots,n$. 
However, \eqref{lyapaverasymp} is not suitable for input-to-state stability (ISS) analysis presented in \cite[Chapter 4.9]{khalil} due to the asymptotically decaying signal in the third term. Nevertheless, the analysis can be performed in a different time domain.
Let us define a new time domain as follows 
\begin{align}
    \tau={}&t_0+\gamma \left(\xi^{-2\kappa+r+v}(t)-1\right), \qquad \tau \in [t_0, \infty), \nonumber \\
    ={}&t_0+\gamma\left((1+\beta (t-t_0))^{\frac{-2\kappa+r+v}{v}}-1\right), 
\end{align}
where $\gamma=\frac{v}{\beta(-2\kappa+r+v)}>0$.
Then, we get
\begin{align}
    \frac{d\tau}{dt}={}(1+\beta (t-t_0))^{\frac{-2\kappa+r}{v}}=\xi^{-2\kappa+r}(t).
\end{align}
We rewrite \eqref{lyapaverasymp} in the contracted time domain $\tau$ as follows
\begin{align}
    \frac{dV}{d\tau}\leq{}&  \xi_{\tau}^{2\kappa-r}(\tau) |\bar{\theta}_f|\left(M_{\theta} \xi_{\tau}^{1+c}(\tau)  +\frac{\beta}{v} \xi_{\tau}^{-v}(\tau) |\bar{\theta}_f|\right) \nonumber \\
    &-\frac{ (k \alpha)_{\text{min}}}{2}\rho_1|\bar{\theta}_f|^{2\kappa}, \label{dvdtauintau}
\end{align}
where
\begin{align}
    \xi_{\tau}(\tau)={}\left(1+(\tau-t_0)/\gamma\right)^{\frac{1}{-2\kappa+r+v}}.
\end{align}
For $\kappa=1$, \eqref{dvdtauintau} results in
\begin{align}
    \frac{dV}{d\tau} \leq{}&-\left(\frac{(k \alpha)_{\text{min}}\rho_1}{2}-\xi_{\tau}^{2-r-v}(\tau)\frac{2\beta}{v}\right)V , \nonumber \\
    &\hspace{0.5cm} \forall |\bar{\theta}_f| \geq 
    \frac{2M_{\theta}\xi_{\tau}^{2-r+1+c}(\tau)}{(k \alpha)_{\text{min}}\rho_1}=:\mathcal{Y}_1\left( \xi_{\tau}^{-1}(\tau) \right), \label{dvdtaukappa1}
\end{align}
from which, by comparison principle, we derive
\begin{align}
    V(\tau)\leq{}&\xi_{\tau}^{\frac{2\beta \gamma }{v}}(\tau)e^{-\frac{ (k \alpha)_{\text{min}}\rho_1}{2}(\tau-t_0)}V(t_0),  \nonumber \\
    &\hspace{2.5cm} \forall |\bar{\theta}_f| \geq 
    \mathcal{Y}_1\left( \xi_{\tau}^{-1}(\tau) \right). \label{Vtaucomp}
\end{align}
Here, $\mathcal{Y}_1$ is a class $\mathcal{K}$ function. For $\kappa>1$, we first apply Young's inequality to get
\begin{align}
    \frac{\beta}{v}|\bar{\theta}_f|\leq{}& \frac{\varsigma  ^{2\kappa-1}}{2\kappa-1}|\bar{\theta}_f|^{2\kappa-1}+M_{\varsigma},
\end{align}
where
\begin{align}
    M_{\varsigma}={}\frac{2\kappa-2}{2\kappa-1}\left(\frac{\beta}{\varsigma v}\right)^{\frac{2\kappa-1}{2\kappa-2}}
\end{align}
with
\begin{align}
    \varsigma={}&\left((2\kappa-1)\frac{ (k \alpha)_{\text{min}}}{8}\rho_1\right)^{\frac{1}{2\kappa-1}},
\end{align}
and rewrite \eqref{dvdtauintau} as
\begin{align}
    \frac{dV}{d\tau}\leq {}&-\frac{ (k \alpha)_{\text{min}}}{4}\rho_1|\bar{\theta}_f|^{2\kappa}-|\bar{\theta}_f|\bigg(\frac{ (k \alpha)_{\text{min}}}{8}\rho_1|\bar{\theta}_f|^{2\kappa-1} \nonumber \\
    &-\xi_{\tau}^{2\kappa-r+1+c}(\tau)M_{\theta}-\xi_{\tau}^{2\kappa-r-v}(\tau)M_{\varsigma}\bigg), \nonumber \\
    \leq{}& -\frac{ (k \alpha)_{\text{min}}}{4}\rho_1|\bar{\theta}_f|^{2\kappa}, \label{dvdtauintaucase2}
    \end{align}
    for all
    \begin{align}
    |\bar{\theta}_f| & \geq  \left(\frac{8 \xi_{\tau}^{2\kappa-r}(\tau)}{(k \alpha)_{\text{min}}\rho_1 }\left(\xi_{\tau}^{1+c}(\tau)M_{\theta}+\xi_{\tau}^{-v}(\tau)M_{\varsigma}\right) \right)^{\frac{1}{2\kappa-1}}\nonumber \\
    &=:\mathcal{Y}_2\left( \xi_{\tau}^{-1}(\tau) \right).
\end{align}
Note that $\mathcal{Y}_2$ is a class $\mathcal{K}$ function. Then, from \eqref{Vtaucomp} and \eqref{dvdtauintaucase2}, we establish the ISS of the $\bar{\theta}_f$-system for $\kappa \geq 1$ in $\tau$-domain. Since $\xi_{\tau}^{2\kappa-r+1+c}(\tau)$ and $\xi_{\tau}^{2\kappa-r-v}(\tau)$ decay to zero due to the conditions $c+1+2\kappa-r<0$ and $2\kappa-r-v<0$, we prove the asymptotic convergence of $\bar{\theta}_f$ to zero in $\tau$ domain, and consequently, in $t$-domain. By exploiting the decaying nature of the inputs of the functions $\mathcal{Y}_1$ and $\mathcal{Y}_2$, we can characterize the ``fading memory'' ISS bound of $\bar{\theta}_f$ (see \cite{grune2002input} and \cite{karafyllis2021input}). Considering \eqref{dvdtaukappa1} and \eqref{dvdtauintaucase2}, the bound in the $t$-domain is given by
\begin{align}
    |\bar{\theta}_f(t)| \leq \mathcal{B}\left(|\bar{\theta}_f(t_0)|, t-t_0\right)+\sup_{t_0 \leq s \leq t}\mathcal{B}\left(\mathcal{Y}\left( \xi^{-1}(s) \right), t-s\right), \label{ISSthetaf}
\end{align}
where $\mathcal{B}$ and $\mathcal{Y}$ are class $\mathcal{KL}$ and $\mathcal{K}$ functions, respectively, defined by
\begin{align}
    &\mathcal{B}\left(|\bar{\theta}_f(t_0)|, t-t_0\right) \nonumber \\
    &={}\begin{cases}
    p_0 e^{-p_1(\xi^{-2\kappa+r+v}(t)-1)}|\bar{\theta}_f(t_0)|,& \text{if } \kappa= 1,\\
    \sqrt{2}(p_2+p_3(\xi^{-2\kappa+r+v}(t)-1))^{\frac{1}{2-2\kappa}},  & \text{if } \kappa>1, 
\end{cases} \label{ISSbound1asympt}
\end{align}
and
\begin{align}
&\mathcal{Y}\left( s \right)={}\begin{cases}
    \mathcal{Y}_1\left( s \right),& \text{if } \kappa= 1,\\
    \mathcal{Y}_2\left( s \right),  & \text{if } \kappa>1, 
\end{cases} \label{ISSbound2asympt}
\end{align}
with some $p_0>0$, $p_1<(k\alpha)_{\text{min}}\rho_1 \gamma/4$, $p_2=\left(|\bar{\theta}_f(t_0)|/2\right)^{1-\kappa}$, $p_3=(k\alpha)_{\text{min}}\rho_1 \gamma (\kappa-1)/4$.

In addition to establishing the asymptotic stability of the $\bar{\theta}_f$-system, we also need to confirm the asymptotic stability of the $\bar{\eta}_f$-system in \eqref{thetaasymlieball}. We first examine the unforced system $\dot{\bar{\eta}}_f={}\left(\frac{2\beta \kappa}{v \xi^v(t)}-\omega_h \right)\bar{\eta}_f$, which yields the following solution 
\begin{align}
    \bar{\eta}_f(t)=\xi^{2\kappa}(t)e^{-\omega_h(t-t_0)} \bar{\eta}_f(t_0).
\end{align}
Therefore, the unforced system is exponentially stable at the origin. Furthermore, by revisiting the bound \eqref{Jfxibound}, we characterize an upper bound for the input of the $\bar{\eta}_f$-system in \eqref{thetaasymlieball} as
\begin{align}
    |\omega_h \xi^{2\kappa}(t) J_f(\bar{\theta}_f,\xi(t)) &-\xi^{2\kappa}(t)\dot{J}(\theta^*(t))| \nonumber \\
    \leq{}& \omega_h \rho_2 |\theta_f|^{2\kappa}+M_J \xi^{d+2\kappa}(t).
\end{align}
Since $d<-2\kappa$ and $|\bar{\theta}_f| \to 0$, the $\bar{\eta}_f$-system is input-to-state stable and uniformly asymptotically stable at the origin.

\textbf{Step 5: Lie bracket averaging theorem.}
Given the uniform asymptotic stability established for the averaged system in \eqref{thetaasymlieball} in Step 4, we conclude from Theorem \ref{LieBracketAvThe}
that the origin of the transformed system \eqref{thetaasymtrans} is practically uniformly asymptotically stable. The existence of $\omega^*$ and the role of $\omega$ are as defined in Definition \ref{def:localpractical}.

\textbf{Step 6: Convergence to extremum.} Considering the
result in Step 5 and recalling from \eqref{xidotdyn}, \eqref{statetransfor} that
\begin{align}
    \theta=\theta^*(t)+\frac{1}{(1+\beta(t-t_0))^{\frac{1}{v}}}\theta_f,
\end{align}
we conclude the asymptotic convergence of $\theta(t)$ to $\theta^*(t)$ at the rate of $1/\xi(t)$. Considering the fading memory ISS bound \eqref{ISSthetaf}, we can provide bound on the convergence error as follows
\begin{align}
    |\theta(t)-\theta^*(t)|\leq{} \xi^{-1}(t)\left(|\theta_f(t)-\bar{\theta}_f(t)|+|\bar{\theta}_f(t)|\right), 
\end{align}from which we conclude \eqref{KLforasymp}. The convergence of $\theta(t)$ to $\theta^*(t)$ implies the asymptotic convergence of the output $y(t)$ and the filtered state $\eta(t)$ to $J(\theta^*(t))$ at the rate of $1/\xi^{2\kappa}(t)$, based on \eqref{Jboundcon}, \eqref{statetransfor}, \eqref{statetransfor2}, and thereby concludes the proof of Theorem \ref{theoremasymp}. \hfill
\end{proof}

\subsection{Asymptotic uES with linearly-chirped probing}
The following theorem presents an asymptotic ES design that employs a chirped probing signal, which grows asymptotically, in contrast to the constant probing in \eqref{ESasymp}.

\begin{theorem}
Consider the following asymptotic uES design
\begin{equation} 
\begin{cases}
\begin{aligned}
\dot{\theta}={}&\xi^{q-v-1}(t) \sum_{i=1}^n \sqrt{\alpha_i \omega_i}  e_i  \\
&\times \cos\Big(\omega_i (t_0+\gamma(\xi^{q}(t)-1))+k_i \xi^r(t)(J(\theta)-\eta)\Big),   \\
\dot{\eta}={}&\left(-\omega_h \eta+\omega_h J(\theta)\right)\xi^{q-v}(t),
\end{aligned}
\end{cases}
\end{equation}
with the asymptotically growing function $\xi(t)$ defined by \eqref{xidotdyn}. Let $\omega_i = \omega \hat{\omega}_i$ such that $\hat{\omega}_i \neq \hat{\omega}_j$ $\forall i \neq j$, $ t_0 \geq 0$, $\alpha_i, k_i, \beta, \omega_h>0$ $\forall i=1,\dots,n$, $\gamma=v/(\beta q)$, and $q > 2\kappa-r \geq 0$. 
Let Assumptions \ref{Ass0}, \ref{Ass1} hold, and Assumption \ref{asympextbound} holds with $\phi(t)=\xi(t)$, $c<q-v-1-2\kappa+r$, $d<q-v-2\kappa$.
Then, there exists $\omega^*>0$
such that for all $\omega > \omega^*$, the input $\theta(t) \to \theta^*(t)$ semi-globally with respect to $\omega$ and asymptotically at the rate of $1/\xi(t)$, and there exist a class $\mathcal{KL}$ function $\mathcal{B}$, a class $\mathcal{K}$ function $\mathcal{Y}$, and a nonnegative constant $D(\theta(t_0), \theta^*(t_0), \eta(t_0))$ such that
    \begin{align}
    |\theta(t)-\theta^*(t)|\leq{}&\xi^{-1}(t) \Big(D+\mathcal{B}\left(|\theta(t_0)-\theta^*(t_0)|, t-t_0\right) \nonumber \\
    &+\sup_{t_0 \leq s \leq t}\mathcal{B}\left(\mathcal{Y}\left( \xi^{-1}(s) \right), t-s\right)\Big). \label{KLforasympchirp}
\end{align}
\end{theorem}

\begin{proof}
Refer to the proof of Theorem \ref{generalproof} by choosing the parameters and satisfying conditions as in Table \ref{table_chirpy}.
\end{proof}

\subsection{Numerical Simulation}
\begin{figure}[t]
    \centering
    \begin{subfigure}[b]{\linewidth}
        \centering
        \includegraphics[width=\linewidth]{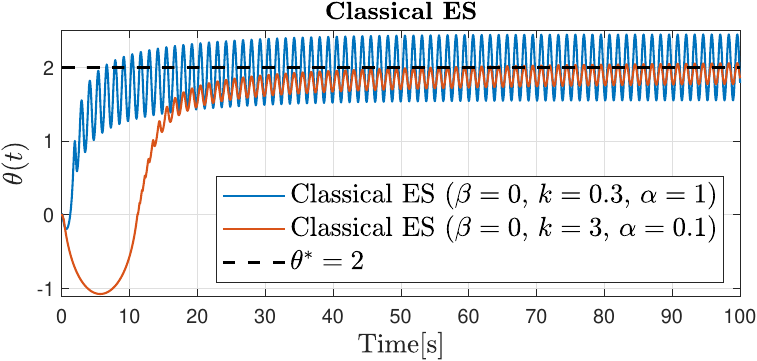}
        \caption{}
        \label{fig:asymp_uES_a}
    \end{subfigure} \\ \vspace{0.2cm}
    \begin{subfigure}[b]{\linewidth}
        \centering
        \includegraphics[width=\linewidth]{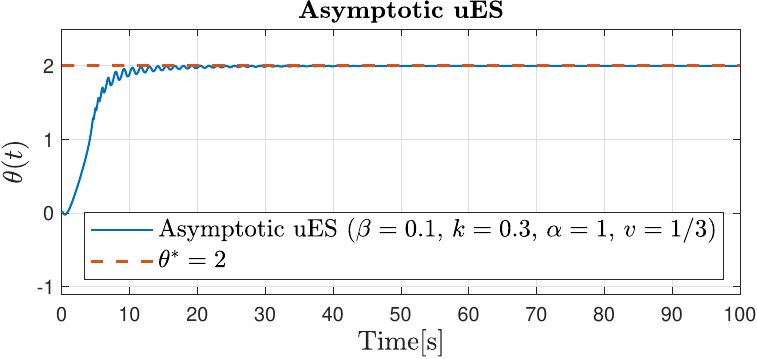}
        \caption{}
        \label{fig:asymp_uES_b}
    \end{subfigure}
    \caption{$(a)$ The trajectories of the classical ES \cite{scheinker2014extremum} with two distinct parameter sets. $(b)$ The trajectory of the asymptotic uES \eqref{ESasymp} using constant frequency probing and achieving a convergence rate of $1/(1+0.1t)^{3}$.}
    \label{fig:asymp_uES}
\end{figure}
This subsection presents numerical results to illustrate the capabilities of the developed designs under two scenarios: a constant optimum and a diverging optimum.
\subsubsection{Seeking constant $\theta^*$} \label{asymuES_numsim1}
We consider the optimization problem \eqref{Jopt} with the following nonlinear map
\begin{align}
    J(\theta)={}1+(\theta-2)^4. \label{funcexamp}
\end{align}
Note that the function \eqref{funcexamp} has a unique minimum at $\theta^*=2$ and satisfies Assumption \ref{Ass1} with $\kappa=2$. 
To provide a basis for comparison, we use the classical ES introduced in \cite{scheinker2014extremum}. Setting $\beta=0$ in the asymptotic uES \eqref{ESasymp}, our design simplifies to the classical approach. In all simulations, initial conditions are set to zero, and the oscillation frequency and high-pass filter frequency are set to $\omega=5$ and $\omega_h=3$, respectively.  

In Figure \ref{fig:asymp_uES_a}, we depict the trajectories of the classical ES with two distinct parameter sets. We observe that the classical ES with $k=0.3$ and $\alpha=1$ converges to a large neighborhood of the optimum due to the high $\alpha$. To reduce the size of the steady-state oscillation, one might consider decreasing $\alpha$ from $1$ to $0.1$, while simultaneously increasing the gain $k$ from $0.3$ to $3$ to maintain the same convergence rate. However, as observed in Figure \ref{fig:asymp_uES_a}, such adjustment leads to poorer transient performance, with an initial deviation in the opposite direction, despite the reduced oscillations at the steady state compared to one with higher amplitude. Our design, illustrated in Figure \ref{fig:asymp_uES_b}, addresses this issue. We employ the asymptotic uES \eqref{ESasymp} with $\beta=0.1$, $k=0.3$, $\alpha=1$, $v=1/3$, $r=4$. It starts with a high $\alpha=1$ and low $k=0.3$, and as $\xi(t)$ increases over time, the input settles to its optimum value with good transient performance. The chosen parameters ensure a convergence rate of $1/(1+0.1t)^{3}$.

\subsubsection{Tracking divergent $\theta^*(t)$} 
\begin{figure}[t]
    \centering
    \includegraphics[width=\columnwidth]{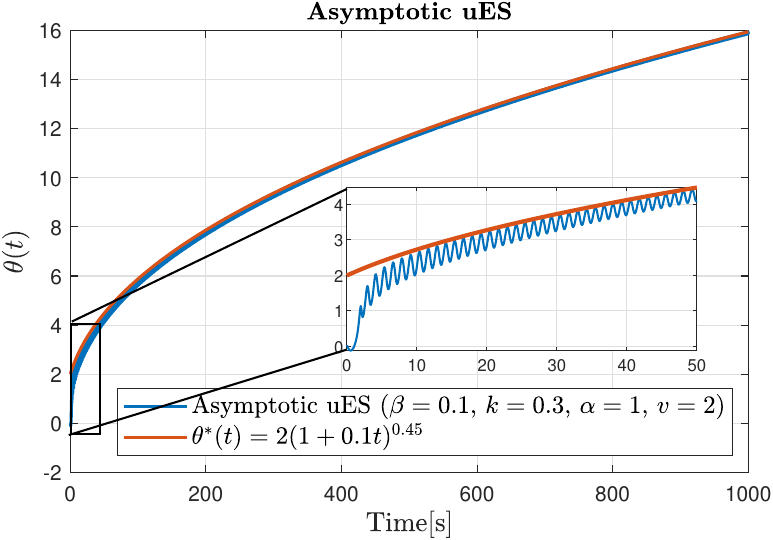} 
    \caption{Perfect tracking of an asymptotically shifting optimum by asymptotic uES with constant-frequency probing \eqref{ESasymp}. The convergence rate is $1/(1+0.1t)^{0.5}$.}
    \label{asymp_uES_track}
\end{figure}
We consider the following map where the optimum input $\theta^*(t)$ shifts from $2$ to infinity
\begin{align}
    J(\theta)={}1+\left(\theta-2(1+0.1t)^{0.45}\right)^4. 
\end{align}
The initial conditions and the parameters are same as in Fig. \ref{fig:asymp_uES}, except that the parameter $v$ is increased from $1/3$ to $v=2$ to make our design \eqref{ESasymp} capable of tracking a slowly shifting optimum. This adjustment slows the convergence rate to $1/(1+0.1t)^{0.5}$. Note that $\theta^*(t)=2(1+0.1t)^{0.45}=2\xi^{\frac{0.9}{2}}(t)$ satisfies \eqref{ass3bound} with $c=-1.1<-1$. Fig. \ref{asymp_uES_track} illustrates the the performance of the design in terms of tracking a varying optimum. 

\section{Exponential uES} \label{expuES}

\subsection{Exponential uES with constant-frequency probing} \label{subsec:expuESconsfreq}
An ES design with unbiased and exponential convergence is presented in the following theorem.

\begin{theorem} \label{theoremexp}
Consider the following exponential uES design
\begin{equation} \label{ESexpoptcase1}
\begin{cases}
\begin{aligned}
\dot{\theta}={}&\zeta^{-1}(t)\sum_{i=1}^n \sqrt{\alpha_i \omega}_i e_i \cos\left(\omega_i t +k_i \zeta^2(t) (J(\theta)-\eta)\right),   \\
\dot{\eta}={}&-\omega_h \eta+\omega_h J(\theta),
\end{aligned}
\end{cases}
\end{equation}
with the exponentially growing function
\begin{align}
    \zeta(t)={}e^{\lambda(t-t_0)}, \quad \lambda>0, \quad t \in [t_0, \infty). \label{xiexp}
\end{align}
Let $\omega_i = \omega \hat{\omega}_i$ such that $\hat{\omega}_i \neq \hat{\omega}_j$ $\forall i \neq j$, $t_0 \geq 0$, and $k_i \alpha_i >{} 2\lambda/\rho_1$, and $\omega_h>{}2\lambda$
for all $i=1,\dots,n$. 
Let Assumptions \ref{Ass0}, \ref{Ass1} holds with $\kappa=1$, and Assumption \ref{asympextbound} holds with $\phi(t)=\zeta(t)$, $c<-1, d<-2$. Then, there exists $\omega^*>0$ such that for all $\omega > \omega^*$, the following holds
\begin{itemize}
    \item $\theta(t) \to \theta^*(t)$ semi-globally with respect to $\omega$ and exponentially at the rate of $\lambda$, and there exist a nonnegative constant $D(\theta(t_0), \theta^*(t_0), \eta(t_0))$ such that 
    \begin{align}
     |\theta(t)-\theta^*(t)|\leq{}&e^{-\lambda(t-t_0)}\bigg[e^{-\ell_0\ell_1(t-t_0)}|\theta(t_0)-\theta^*(t_0)|\nonumber \\
     &+\frac{M_{\theta}}{(1-\ell_0) \ell_1}e^{-\min\{\ell_0\ell_1,\ell_2\}(t-t_0)}+D \bigg] \label{KLforexp}
\end{align}
for any $0<\ell_0<1$ with
\begin{align}
    \ell_1=\,{}&(k\alpha)_{\text{min}}\rho_1/2-\lambda>0, \label{l1def} \\
    \ell_2=\,{}&-\lambda(c+1)>0.\label{l2def}
\end{align}
\item $\eta(t), y(t) \to J(\theta^*(t))$ semi-globally with respect to $\omega$ and exponentially  at the rate of $2\lambda$.
\end{itemize}
\end{theorem}
\begin{remark}
The design parameters $\lambda$ and $(k \alpha)_{\text{min}}$, as well as the system parameters $M_{\theta}$, $\rho_1$, and $c$, determine the size and decay rate of the bound \eqref{KLforexp}. Increasing $\lambda$ accelerates the decay of the update gain and the growth of the controller gain in \eqref{ESexpoptcase1}, leading to a faster convergence rate for the bound \eqref{KLforexp}. The convergence rate $\lambda$ can be made arbitrarily fast, provided that $\lambda < (k \alpha)_{\text{min}} \rho_1 / 2$. A large gain $(k \alpha)_{\text{min}}$ further contributes to a faster decay of the initial conditions and to the rate of  expiration of the time-dependent  term parameterized by $M_{\theta}$ if $\ell_0 \ell_1 < \ell_2$, where $\ell_2$ is an increasing function of $c$. The effect of the residual term $D$, which arises from the difference between the average and original state, diminishes exponentially at the rate of $\lambda$, independent of other parameters. 
\end{remark}
\begin{proof}
Let us proceed through the proof step by step. 

\textbf{Step 1: State transformation.} Let us consider the following transformations
\begin{align}
    \theta_{f}={}& \zeta(t) (\theta-\theta^*(t)), \label{statetransforexp} \\
    \eta_{f}={}& \zeta^{2}(t) (\eta-J(\theta^*(t))). \label{filttransforexp}
\end{align}
Using \eqref{statetransforexp} and \eqref{filttransforexp}, we transform \eqref{ESexpoptcase1} to
\begin{equation} \label{exptransfcas1}
\begin{cases}
\begin{aligned}
\dot{\theta}_{f}={}&-\zeta(t)\dot{\theta}^*(t)+\lambda {\theta}_{f}+\sum_{i=1}^n  \sqrt{\alpha_i \omega_i} e_i \nonumber \\
&\times \cos\Big(\omega_i t+k_i \zeta^2(t)  \left(J_f(\theta_f, t)-\zeta^{-2}(t){\eta}_f\right)\Big), \\
    \dot{\eta}_{f}={}&(2\lambda -\omega_h )\eta_f+\omega_h \zeta^{2}(t)J_f(\theta_f, t)-\zeta^2(t)\dot{J}(\theta^*(t)),
\end{aligned}
\end{cases}
\end{equation}
with
\begin{align}
    J_f(\theta_f, t) ={}& J(\theta_f/\zeta(t) +\theta^*(t))-J(\theta^*(t)). \label{Jfzetat}
\end{align}

\textbf{Step 2: Lie bracket averaging and stability analysis.} 
The feasibility of the error system \eqref{exptransfcas1} for Lie bracket averaging can be verified analogously to Step 2 in the proof of Theorem \ref{theoremasymp}. In contrast to \eqref{thetaasymlieball}, the average of \eqref{exptransfcas1} results in
\begin{equation} \label{Lieexpkap1}
\begin{cases}
\begin{aligned}
\dot{\bar{\theta}}_{f}={}&-\zeta(t)\dot{\theta}^*(t)+\lambda {\bar{\theta}}_{f}-\sum_{i=1}^n \frac{k_i \alpha_i}{2} e_i \zeta^{2}(t) \frac{\partial J_f(\bar{\theta}_f, t)}{\partial \bar{\theta}_{f,i}}, \\
\dot{\bar{\eta}}_{f}={}&(2\lambda -\omega_h )\bar{\eta}_f+\omega_h \zeta^{2}(t)J_f(\bar{\theta}_f, t)-\zeta^2(t)\dot{J}(\theta^*(t)).
\end{aligned}
\end{cases}
\end{equation}
Consider the following Lyapunov function
\begin{align}
    V=\frac{1}{2} |\bar{\theta}_f|^2. \label{Lyapkeq1}
\end{align}
Taking Assumption \ref{Ass1} into consideration and using \eqref{Lieexpkap1}, we compute the time derivative of \eqref{Lyapkeq1} as
\begin{align}
    \dot{V} \leq{}&\zeta(t)|\bar{\theta}_f| |\dot{\theta}^*(t)|+\lambda |\bar{\theta}_f|^2-\sum_{i=1}^n \frac{k_i \alpha_i}{2}  \zeta^{2}(t)  \bar{\theta}_{f,i}\frac{\partial J_f(\bar{\theta}_f, t)}{\partial \bar{\theta}_{f,i}}, \nonumber \\
    \leq{} & M_{\theta} \zeta^{c+1}(t) |\bar{\theta}_f| -
    \left(\frac{ (k \alpha)_{\text{min}}\rho_1}{2}-\lambda \right)|\bar{\theta}_f|^2,
    \nonumber \\
    \leq{} & -\ell_0\ell_1 |\bar{\theta}_f|^2, \qquad \forall |\bar{\theta}_f| \geq \frac{M_{\theta}\zeta^{c+1}(t)}{(1-\ell_0)\ell_1}=:\mathcal{Y}\left( \zeta^{-1}(t)\right), \label{VdotISSexp}
\end{align}
which establishes ISS of the $\bar{\theta}_f$-system for $k_i \alpha_i>2\lambda/\rho_1$. In \eqref{VdotISSexp}, $0<\ell_0<1$, $\ell_1$ is defined in \eqref{l1def}, and $\mathcal{Y}$ is a class $\mathcal{K}$ function. Recalling $c<-1$, we prove the asymptotic convergence of $\bar{\theta}_f$ to zero. The corresponding fading memory ISS bound can be written as 
\begin{align}
    |\bar{\theta}_f(t)| \leq \mathcal{B}\left(|\bar{\theta}_f(t_0)|, t-t_0\right)+\sup_{t_0 \leq s \leq t}\mathcal{B}\left(\mathcal{Y}\left( \zeta^{-1}(s) \right), t-s\right) \label{ISSthetafexp}
\end{align}
with
\begin{align}
\mathcal{B}\left(|\bar{\theta}_f(t_0)|, t-t_0\right)={}&
    e^{-\ell_0\ell_1(t-t_0)}|\bar{\theta}_f(t_0)|. \label{bdefforexp}
\end{align}
We can characterize the decay rate of the sup term in \eqref{ISSthetafexp} using \eqref{xiexp}, \eqref{VdotISSexp} and \eqref{bdefforexp} as follows
\begin{align}
    \hspace{-0.14cm}\sup_{t_0 \leq s \leq t}\mathcal{B}\left(\mathcal{Y}\left( \zeta^{-1}(s) \right), t-s\right)=&\sup_{t_0 \leq s \leq t} \frac{M_{\theta}e^{-\ell_0\ell_1(t-s)-\ell_2 (s-t_0)}}{(1-\ell_0)\ell_1} \nonumber \\
     \leq &\, \frac{M_{\theta}e^{-\min\{\ell_0\ell_1,\ell_2\}(t-t_0)}}{(1-\ell_0)\ell_1}
\end{align}
with $\ell_2$ as defined in \eqref{l2def}, respectively. We note that $-\ell_0\ell_1(t-s)-\ell_2 (s-t_0)<-\ell_0\ell_1(t-t_0)$ for $\ell_0\ell_1 \leq \ell_2$ and $-\ell_0\ell_1(t-s)-\ell_2 (s-t_0)<-\ell_2(t-t_0)$ for $\ell_2 \leq \ell_0\ell_1$.

To establish the asymptotic stability of the $\bar{\eta}_f$-system, 
we analyze its unforced dynamics, governed by $\dot{\bar{\eta}}_f=(2\lambda-\omega_h)\bar{\eta}_f$. This system is exponentially stable at the origin for $\omega_h>2\lambda$.
The input-to-state stability of the $\bar{\eta}_f$-system is established by noting from \eqref{Jboundcon} that
\begin{align}
     \zeta^{2}(t)|\omega_h J_f(\bar{\theta}_f,\zeta(t))-\dot{J}(\theta^*(t))| 
    \leq{}& \omega_h \rho_2 |\bar{\theta}_f|^2 +M_{\theta} \zeta^{d+2}(t).   
\end{align}
Since $d<-2$, the asymptotic convergence of $\bar{\theta}_f$ results in the asymptotic convergence of $\bar{\eta}_f$ as well.

\textbf{Step 3: Lie bracket averaging theorem.}
Given the uniform asymptotic stability established for the averaged system in \eqref{Lieexpkap1} in Step 2, we conclude from Theorem \ref{LieBracketAvThe}
that the origin of the transformed system \eqref{exptransfcas1} is practically uniformly asymptotically stable.
The existence of $\omega^*$ and the role of $\omega$ are as defined in Definition \ref{def:localpractical}.

\textbf{Step 4: Convergence to extremum.} Considering the result in Step 3 and recalling from \eqref{xiexp}, \eqref{statetransforexp} that
\begin{align}
    \theta=\theta^*(t)+ e^{-\lambda(t-t_0)}\theta_f,
\end{align}
we conclude the exponential convergence of $\theta(t)$ to $\theta^*(t)$ at the rate of $\lambda$. Recalling the ISS bound \eqref{ISSthetafexp}, the bound on the convergence error is provided below
\begin{align}
    &|\theta(t)-\theta^*(t)| \leq{} \zeta^{-1}(t)\left(|\theta_f(t)-\bar{\theta}_f(t)|+|\bar{\theta}_f(t)|\right), 
\end{align}from which we conclude \eqref{KLforexp}.
The convergence of $\theta(t)$ to $\theta^*(t)$ implies the exponential convergence of the output $y(t)$ and the filtered state $\eta(t)$ to $J(\theta^*(t))$ at the rate of $2\lambda$, as evident from \eqref{Jboundcon}, \eqref{statetransforexp}, \eqref{filttransforexp},
and completes the proof of Theorem \ref{theoremexp}.
\hfill
\end{proof}

\subsection{Exponential uES with exponentially-chirped probing} \label{subsec:expuESchirped}
The following theorem presents an alternative design to \eqref{ESexpoptcase1} using an exponentially chirped probing signal.

\begin{theorem}
Consider the following asymptotic uES design
\begin{equation} \label{expuESchirpy}
\begin{cases}
\begin{aligned} 
\dot{\theta}={}&\zeta^{q-1}(t) \sum_{i=1}^n \sqrt{\alpha_i \omega_i}  e_i  \\
&\times \cos\Big(\omega_i (t_0+\gamma(\zeta^{q}(t)-1))+k_i \zeta^r(t)(J(\theta)-\eta)\Big),   \\
\dot{\eta}={}&\left(-\omega_h \eta+\omega_h J(\theta)\right)\zeta^{q}(t),
\end{aligned}
\end{cases}
\end{equation}
with the exponentially growing function $\zeta(t)$ defined by \eqref{xiexp}. Let $\omega_i = \omega \hat{\omega}_i$ such that $\hat{\omega}_i \neq \hat{\omega}_j$ $\forall i \neq j$, $ t_0 \geq 0$, $\alpha_i, k_i, \beta, \omega_h>0$ $\forall i=1,\dots,n$,  $\gamma=1/(\lambda q)$, and $q > 2\kappa-r \geq 0$.  
Let Assumptions \ref{Ass0}, \ref{Ass1} hold, and Assumption \ref{asympextbound} holds with $\phi(t)=\zeta(t)$, $c<q-1-2\kappa+r$, $d<q-2\kappa$.
Then, there exists $\omega^*>0$
such that for all $\omega > \omega^*$, the input $\theta(t) \to \theta^*(t)$ semi-globally with respect to $\omega$ and exponentially at the rate of $\lambda$, and there exist a class $\mathcal{KL}$ function $\mathcal{B}$, a class $\mathcal{K}$ function $\mathcal{Y}$, and a nonnegative constant $D(\theta(t_0), \theta^*(t_0), \eta(t_0))$ such that 
    \begin{align}
    |\theta(t)-\theta^*(t)|&\leq{}\zeta^{-1}(t) \Big(D+\mathcal{B}\left(|\theta(t_0)-\theta^*(t_0)|, t-t_0\right) \nonumber \\
    &\hspace{0.3cm}+\sup_{t_0 \leq s \leq t}\mathcal{B}\left(\mathcal{Y}\left( \zeta^{-1}(s) \right), t-s\right)\Big). \label{KLforexpchirp}
\end{align}
\end{theorem}

\begin{proof}
Refer to the proof of Theorem \ref{generalproof} by choosing the parameters and satisfying conditions as in Table \ref{table_chirpy}.
\end{proof}

\subsection{Numerical Simulation}

\begin{figure}[t]
    \centering
    \begin{subfigure}[b]{\linewidth}
        \centering
        \includegraphics[width=\linewidth]{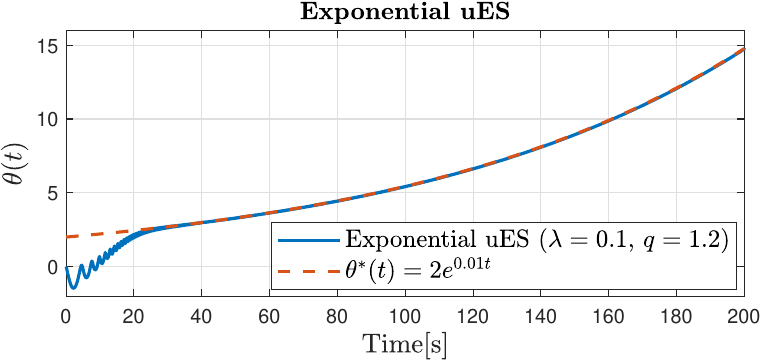}
        \caption{}
        \label{fig:exponent_uES_track_a}
    \end{subfigure} \\ \vspace{0.2cm}
    \begin{subfigure}[b]{\linewidth}
        \centering
        \includegraphics[width=\linewidth]{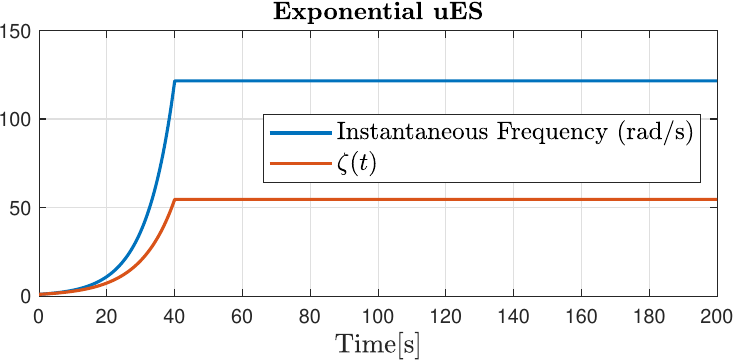}
        \caption{}
        \label{fig:exponent_uES_track_b}
    \end{subfigure}
    \caption{$(a)$ Tracking of an exponentially shifting optimum by exponential uES with chirped probing \eqref{expuESchirpy}. $(b)$ Evolution of $\zeta(t)$ and the instantaneous frequency $\omega \gamma d\zeta^q/dt$, with their growth terminated at $t = 40$ seconds.}
    \label{fig:exponent_uES_track}
\end{figure}

We consider a quadratic map of the form 
\begin{align}
    J(\theta)={}1+\left(\theta-2e^{0.01t}\right)^2,
\end{align}
in which the optimum input $\theta^*(t)=2e^{0.01t}$ diverges exponentially. We apply the exponential uES with chirped probing \eqref{expuESchirpy} by setting the parameters $\lambda=0.1$, $k=0.3$, $\omega=1$, $\omega_h=3$, $\alpha=1$, $q=1.2$, and $r=2$. Note that $\theta^*(t)=2\zeta^{0.1}(t)$ and $\kappa=1$. Thus, the parameters $c$ and $q$ satisfy the conditions $q=1.2>2\kappa-r=0$ and $c=0.1<q-1-2\kappa+r=0.2$. We present the results in Fig. \ref{fig:exponent_uES_track}. As discussed in Remark \ref{remark:saturation}, we terminate the growth of the signal $\zeta(t)$ and the instantaneous frequency, $\omega \gamma d\zeta^q/dt$, at $t=40$ seconds, which is sufficient to track the optimum closely.

\section{Prescribed-Time uES} \label{subsec:prescuES}
We present the prescribed-time uES design in the following theorem, which offers the strongest result in terms of convergence speed and convergence error compared to the asymptotic and exponential designs.

\begin{theorem}
Consider the following prescribed-time uES design
\begin{equation} \label{ptesdesign}
\begin{cases}
\begin{aligned}
\dot{\theta}={}&\mu^{q+\varrho-1}(t) \sum_{i=1}^n \sqrt{\alpha_i \omega_i}  e_i  \\
&\times \cos\Big(\omega_i (t_0+\gamma(\mu^{q}(t)-1))+k_i \zeta^r(t)(J(\theta)-\eta)\Big),   \\
\dot{\eta}={}&\left(-\omega_h \eta+\omega_h J(\theta)\right)\mu^{q+\varrho}(t),
\end{aligned}
\end{cases}
\end{equation}
with the blow-up function
\begin{align}
    \mu(t)={}\left(\frac{T}{T+t_0-t}\right)^{\frac{1}{\varrho}}, \qquad t \in [t_0, t_0+T), \label{mudefinition}
\end{align}
Let $\omega_i = \omega \hat{\omega}_i$ such that $\hat{\omega}_i \neq \hat{\omega}_j$ $\forall i \neq j$, $ t_0 \geq 0$, $\alpha_i, k_i, \beta, \omega_h>0$ $\forall i=1,\dots,n$, $\gamma=\varrho T/q$, and $q > 2\kappa-r \geq 0$.  Let Assumptions \ref{Ass0}, \ref{Ass1} hold, and Assumption \ref{asympextbound} holds with $\phi(t)=\mu(t)$, $c<q+\varrho-1-2\kappa+r$, $d<q+\varrho-2\kappa$.
Then, there exists $\omega^*>0$ such that for all $\omega > \omega^*$, the input $\theta(t) \to \theta^*(t)$ semi-globally with respect to $\omega$ and in prescribed-time $t_0+T$, and there exist a class $\mathcal{KL}$ function $\mathcal{B}$, a class $\mathcal{K}$ function $\mathcal{Y}$, and a nonnegative constant $D(\theta(t_0), \theta^*(t_0), \eta(t_0))$ such that 
    \begin{align}
    |\theta(t)-\theta^*(t)|&\leq{}\mu^{-1}(t) \Big(D+\mathcal{B}\left(|\theta(t_0)-\theta^*(t_0)|, t-t_0\right) \nonumber \\
    &\hspace{0.3cm}+\sup_{t_0 \leq s \leq t}\mathcal{B}\left(\mathcal{Y}\left( \mu^{-1}(s) \right), t-s\right)\Big). \label{KLforpt}
\end{align}
\end{theorem}

\begin{proof}
Refer to the proof of Theorem \ref{generalproof} by choosing the parameters and satisfying conditions as in Table \ref{table_chirpy}.
\end{proof}

\begin{remark}
For strongly convex maps with $\kappa=1$ and fixed optima, two alternative PT-uES designs are developed in \cite{yilmaz2024unbiased}. One design has a frequency that grows at a power-of-$\mu$ rate and an update rate that decays at an exp-of-power-of-$\mu$ rate. The other design has a frequency that grows at a log-of-$\mu$ rate and an update rate that decays at a power-of-$\mu$ rate. However, convex maps with $\kappa \geq 1$ and time-varying optima require a different design approach, as presented in \eqref{ptesdesign}.
\end{remark}

\subsection{Numerical Simulation}
\begin{figure}[t]
    \centering
    \begin{subfigure}[b]{\linewidth}
        \centering
        \includegraphics[width=\linewidth]{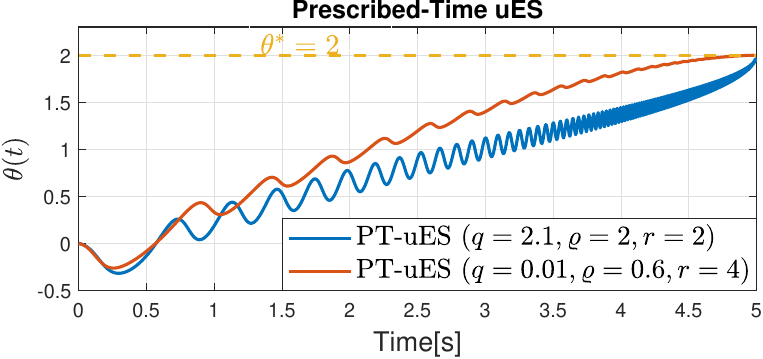}
        \caption{}
        \label{fig:presc_uES_a}
    \end{subfigure} \\ \vspace{0.2cm}
    \begin{subfigure}[b]{\linewidth}
        \centering
        \includegraphics[width=\linewidth]{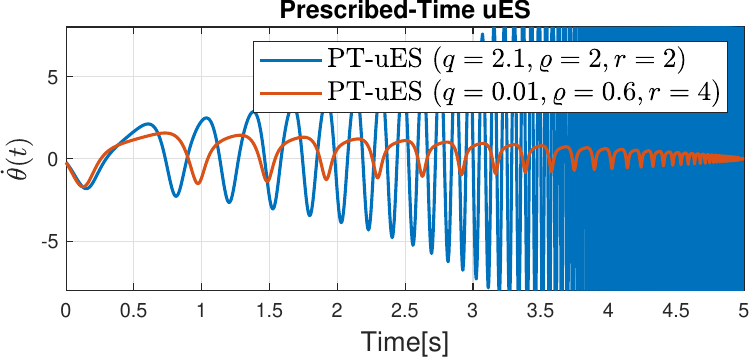}
        \caption{}
        \label{fig:presc_uES_b}
    \end{subfigure}
    \caption{$(a)$ Unbiased convergence to $\theta^*=2$ in user-prescribed $T=5$ seconds with two different parameter sets. $(b)$ Evolution of their corresponding update laws.}
    \label{fig:presc_uES}
\end{figure}
In this section, we return to the cost function \eqref{funcexamp}, now aiming for convergence to $\theta^*$ within the time period of $T=5$ seconds. We apply the prescribed-time uES \eqref{ptesdesign} with two different parameter sets:
\begin{itemize}
    \item In the first configuration, we use parameters $q=0.01$, $\varrho=0.6$ $\alpha_1=1$, $\omega_1=10$, $\omega_h=3$, and  $k_1=0.1$. The motivation here is to demonstrate that the update law, $\dot{\theta}$, can be guaranteed to remain bounded by choosing $r=2\kappa=4$ and making $p=q+\varrho-1$ negative through small values of $q$ and $\varrho$.
    \item In the second configuration, we use parameters $q=2.1$, $\varrho=2$, and the parameters $\alpha_1, \omega_1, \omega_h$, and $k_1$ are the same as in the first configuration. Opting for a ``safe'' $r=2$ (as discussed in Remark \ref{remark:rv}) leads to a compromise on the boundedness of the update law. This is because the condition $q > 2\kappa-r=2$ must be met, which forces $p=q+\varrho-1$ to be positive. To mitigate the aggressive increase in frequency during the transient, a relatively high value of $\varrho$ is used compared to the first configuration.
\end{itemize}

We present the simulation results in Fig. \ref{fig:presc_uES}, which illustrate the trade-off between the parameter sets. Both designs achieve prescribed-time convergence, but they differ in the growth rate of adaptation and frequency.

\section{Unbiased Source Seeking by Unicycle} \label{sec:experiment}

In this section, we investigate the problem of source localization by a unicycle in a two-dimensional plane modeled by
\begin{align}
    \dot{x}={}&u_1 \begin{bmatrix}
        \cos(\theta) \\ \sin(\theta)
    \end{bmatrix}, \label{unicy_xdot} \\
    \dot{\theta}={}&u_2, \label{unicy_thetadot}
\end{align}
where $x=\begin{bmatrix}
    x_1, & x_2
\end{bmatrix}^T \in \mathbb{R}^2$ denotes the coordinates of the vehicle's center with $x(0) = x_0$, $\theta \in \mathbb{R}$ is the orientation of the unicycle with $\theta(0)=\theta_0$, $u_1(t) \in \mathbb{R}$ and $u_2(t) \in \mathbb{R}$ are the forward and angular velocity inputs, respectively. 
For simplicity, but without loss of generality, we assume
that $J(\cdot)$ is quadratic with diagonal Hessian
matrix, i.e.,
\begin{align}
    J(x)={}J^*-\frac{\rho_{x_1}}{2}(x_1-x_1^*)^2-\frac{\rho_{x_2}}{2}(x_2-x_2^*)^2, \label{cost_unic}
\end{align}
where $x^*= \begin{bmatrix}
    x_1^*, & x_2^*
\end{bmatrix}^T \in \mathbb{R}^2$ is the unknown maximizer, $J^*=J(x^*) \in \mathbb{R}$ is the unknown maximum, and $\rho_{x_1}$, $\rho_{x_2}$ are some unknown positive constants. Our aim is to design an unbiased source seeker that drives the vehicle to the exact position of the source, i.e. $x \to x^*$, exponentially at a user-defined rate of $\lambda>0$. 

The source seeker design follows similar steps to those of exponential uES with constant-frequency probing in Section \ref{subsec:expuESconsfreq}, but a challenge arises because the angle $\theta$ is not directly controlled; instead, the angular velocity $\dot{\theta}$ is manipulated. This makes the system \eqref{unicy_xdot} and \eqref{unicy_thetadot} compatible with \eqref{singpert1} and \eqref{singpert2} instead of \eqref{conaff}, requiring singular perturbation analysis as another layer of the design.

\subsection{Design and Analysis}
We design the inputs $u_1$ and $u_2$ in \eqref{unicy_xdot} and \eqref{unicy_thetadot} such that
\begin{align}
    \dot{x}={}&e^{-\lambda t} \sqrt{\alpha \omega}  \begin{bmatrix}
        \cos(\theta) \\ \sin(\theta)
    \end{bmatrix}, \label{xdot_unic} \\
    \dot{\theta}={}&\omega-ke^{2\lambda t}\left(\eta+J(x)\right) \label{thetadot_unic}
\end{align}
with a low-pass filter
\begin{align}
    \dot{\eta}={}-\omega_h\eta-\omega_h J(x). \label{eta_unic}
\end{align}
The parameters satisfy 
\begin{align}
    k \alpha >{}&\frac{2\lambda }{\rho_{x_i}}(\omega_h-2\lambda), \quad \lambda<\frac{\omega_h}{2}, \quad i=1,2. \label{param_unic}
\end{align}
We present the exponential convergence result for this source seeking problem as follows.
\begin{theorem} 
Consider the closed-loop system \eqref{xdot_unic}--\eqref{thetadot_unic} with the parameters that satisfy \eqref{param_unic} and with the cost function of the form \eqref{cost_unic}. Then, there exists $\omega^*$ such that for all $\omega > \omega^*$ there exists $\omega^*_h>0$ such that for all $\omega_h>\omega_h^*+2\lambda$, $x(t)$ semi-globally exponentially converges to $x^*$ at the rate of $\lambda$.
\end{theorem}

\begin{proof}
\textbf{Step 1: State transformation.}
Let us consider the following transformations
\begin{align}
    x_{f}={}& e^{\lambda t} (x-x^*),  \\
    \eta_{f}={}& e^{2\lambda t} \left(\eta+ J^*\right) \label{unicetaf}
\end{align}
with
\begin{align}
    x_f={}&\begin{bmatrix}
        x_{1,f} & x_{2,f}
    \end{bmatrix}^T.
\end{align}
We write the transformed system using \eqref{xdot_unic} and \eqref{eta_unic} as
\begin{align}
    \dot{x}_f={}&\begin{bmatrix}
        \lambda x_{1,f}+\sqrt{\alpha \omega}\cos(\theta) \\
        \lambda x_{2,f}+\sqrt{\alpha \omega}\sin(\theta)
    \end{bmatrix}, \label{xf_unic} \\
    \dot{\theta}={}&\omega-k({\eta}_f+J_f(x_f)), \label{thetaf_unic} \\
    \dot{\eta}_f={}&-(\omega_h-2\lambda)\eta_f-\omega_h J_f(x_{f}), \label{etadf_unic}
\end{align}
where
\begin{align}
    J_f(x_{f})={}&-\frac{\rho_{x_1}}{2}x_{1,f}^2-\frac{\rho_{x_2}}{2}x_{2,f}^2.
\end{align}
In view of \eqref{etadf_unic}, $\theta$-dynamics in \eqref{thetaf_unic} can be expressed as
\begin{align}
    \dot{\theta}={}\omega+\frac{k}{\omega_h}\dot{\eta}_f-\frac{2 k \lambda}{\omega_h}\eta_f. \label{thetadot_unicy}
\end{align}
Considering \eqref{thetadot_unicy}, let us perform the change of variables
\begin{align}
    \theta_e={}\theta-\omega t-\frac{k}{\omega_h}\eta_f,
\end{align}
and rewrite \eqref{xf_unic} as 
\begin{align}
   \dot{x}_f ={}&\begin{bmatrix}
        \lambda x_{1,f}+\sqrt{\alpha \omega}\cos(\omega t+\check{k} \eta_f+\theta_e) \\ \lambda x_{2,f}+\sqrt{\alpha \omega}\sin(\omega t+\check{k} \eta_f+\theta_e) 
    \end{bmatrix}, \label{xfcov} \\
    \dot{\theta}_e={}&-2\check{k}\lambda \eta_f, \\
    \dot{\eta}_f={}&-(\omega_h-2\lambda)\eta_f-(\omega_h-2\lambda) \check{J}_f(x_{f}), \label{etafcov}
\end{align}
where
\begin{align}  
    \check{k}={}&k/\omega_h, \label{kcond} \\
    \check{J}_f(x_f)={}&-\frac{\check{\rho}_{x_1}}{2}x_{1,f}^2-\frac{\check{\rho}_{x_2}}{2}x_{2,f}^2, 
\end{align}
with
\begin{align}
 \check{\rho}_{x_i}={}&\rho_{x_i}\frac{\omega_h}{\omega_h-2\lambda}, \quad i=1,2. \label{rhoxdef}
\end{align}

\textbf{Step 2: Singular perturbation analysis.} In order to use the singular perturbation analysis outlined in Section \ref{tit:sing_pert}, we consider a new time scale $t=\varepsilon \tau$, where $\varepsilon=1/(\omega_h-2\lambda)$. Then, using trigonometric identities, we rewrite \eqref{xfcov}--\eqref{etafcov} in $\tau$-domain as
\begin{align}
    \frac{d}{d \tau}\begin{bmatrix}
        x_{1,f} \\ x_{2,f} \\ \theta_e
    \end{bmatrix} ={}&\varepsilon\begin{bmatrix}
        \lambda {x}_{1,f} \\  \lambda {x}_{2,f} \\ -2\check{k} \lambda \eta_f
    \end{bmatrix}+\varepsilon\begin{bmatrix}
        \sqrt{\alpha \omega}\cos(\check{k}\eta_f+\theta_e) \\  \sqrt{\alpha \omega}\sin(\check{k}\eta_f+\theta_e) \\ 0
    \end{bmatrix} \nonumber \\
    &\times \cos(\omega \varepsilon \tau)+\varepsilon\begin{bmatrix}
        -\sqrt{\alpha \omega}\sin(\check{k}\eta_f+\theta_e) \\  \sqrt{\alpha \omega}\cos(\check{k}\eta_f+\theta_e) \\ 0
    \end{bmatrix}\nonumber \\
    &\times \sin(\omega \varepsilon \tau), \label{xfsingular} \\
    \frac{d{\eta}_f}{d\tau}={}&\underbrace{-\eta_f- \check{J}_f(x_{f})}_{=g(x_f,\eta_f)}. \label{etafsubg}
\end{align}
The quasi-steady-state is $\eta_f=-\check{J}_f(x_{f})$. Defining ${\eta}_{f,b}=\eta_f-(-\check{J}_f(x_f))$, the boundary layer model for \eqref{etafsubg} is obtained es 
\begin{align}
    \frac{d{\eta}_{f,b}}{d\tau}={}&g(x_f,\eta_{f,b}-\check{J}_f(x_f)) \nonumber \\
    ={}&-{\eta}_{f,b}, \label{boundlayer}
\end{align}
which is globally exponentially stable. By substituting
the quasi-steady state into \eqref{xfsingular}, we obtain the reduced model as
\begin{align}
      \frac{d}{d \tau}\begin{bmatrix}
         x_{1,f,r} \\ x_{2,f,r} \\ \theta_{e,r}
     \end{bmatrix}&={}\varepsilon\begin{bmatrix}
         \lambda {x}_{1,f,r} \\  \lambda {x}_{2,f,r} \\ 2\check{k} \lambda \check{J}_f(x_{f,r})
    \end{bmatrix}\nonumber \\
    &+\varepsilon\begin{bmatrix}
        \sqrt{\alpha }\cos(-\check{k}\check{J}_f(x_{f,r})+\theta_{e,r}) \\  \sqrt{\alpha }\sin(-\check{k}\check{J}_f(x_{f,r})+\theta_{e,r}) \\ 0
    \end{bmatrix} \sqrt{\omega} \cos(\omega \varepsilon \tau) \nonumber \\
    &+\varepsilon\begin{bmatrix}
        -\sqrt{\alpha }\sin(-\check{k}\check{J}_f(x_{f,r})+\theta_{e,r}) \\  \sqrt{\alpha }\cos(-\check{k}\check{J}_f(x_{f,r})+\theta_{e,r}) \\ 0
    \end{bmatrix} \sqrt{\omega} \sin(\omega \varepsilon \tau) \label{xfexpanded}
\end{align}
with $x_{f,r}=\begin{bmatrix}
    x_{1,f,r} & x_{2,f,r}
\end{bmatrix}^T$.

\begin{figure*}[ht]
    \centering
    \begin{subfigure}[b]{.3\linewidth}
        \centering
        \includegraphics[trim=0cm 1cm 1.8cm 1cm, clip=true, width=\linewidth]{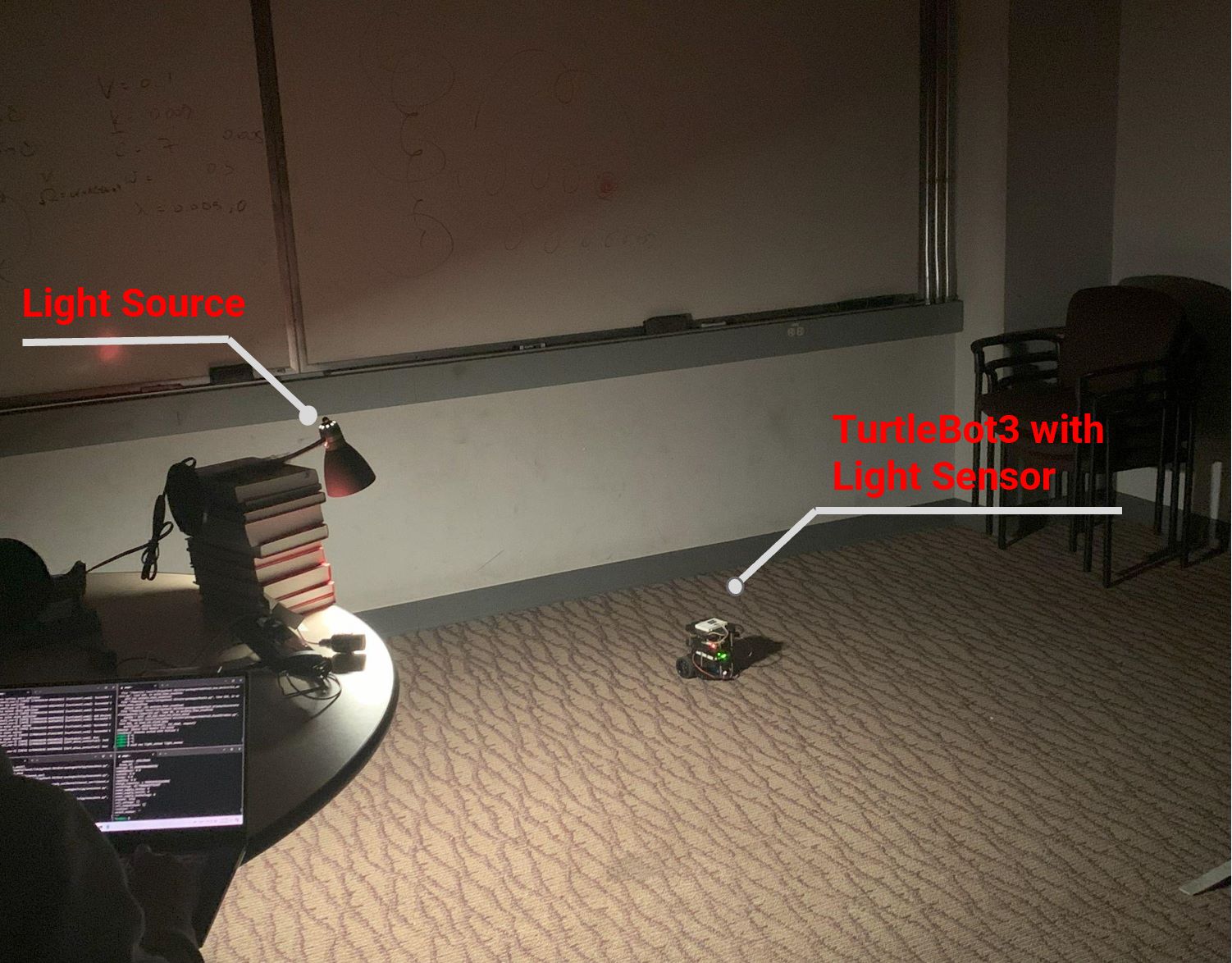} 
        \caption{}
        \label{fig:experiment_sub1}
    \end{subfigure}%
    ~ 
    \begin{subfigure}[b]{.69\linewidth}
        \centering
        \includegraphics[width=\linewidth]{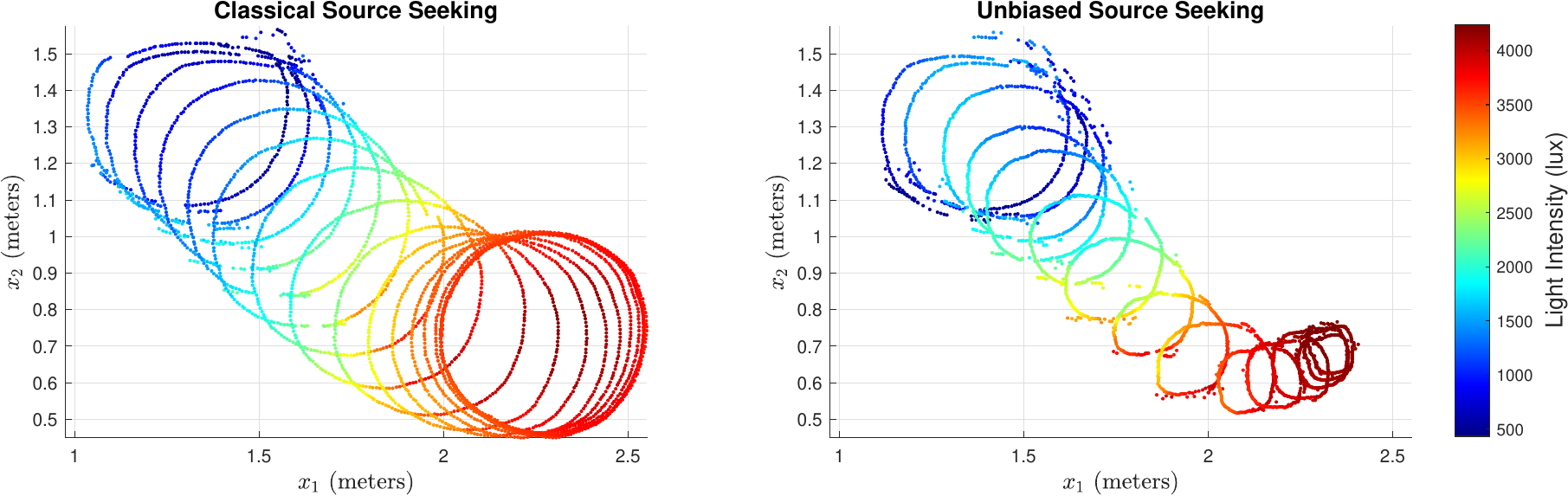}
        \caption{}
        \label{fig:experiment_sub2}
    \end{subfigure}
    \caption{$(a)$ Experimental setup consisting of a light source and a unicycle robot equipped with a light sensor on top. $(b)$ The vehicle trajectory for two ES designs on a 2D plane, with color-coded light intensity in lux representing the intensity at each position.}
    \label{fig:experiment}
\end{figure*}

\textbf{Step 3: Lie bracket averaging and stability analysis.} The Lie bracket average of the reduced system \eqref{xfexpanded} is given by
\begin{align}
     \frac{d}{d \tau}\begin{bmatrix}
         \bar{x}_{1,f,r} \\ \bar{x}_{2,f,r} \\ \bar{\theta}_{e,r}
     \end{bmatrix} ={}&\begin{bmatrix}
         (\lambda-\check{k}  \alpha\check{\rho}_{x_1}/2)\bar{x}_{1,f,r}   \\ (\lambda -\check{k}  \alpha\check{\rho}_{x_2}/2)\bar{x}_{2,f,r} \\ 2\check{k} \lambda \check{J}_f(\bar{x}_{f,r})
     \end{bmatrix}, \label{lieunicy}
\end{align}
with $\bar{x}_{f,r}=\begin{bmatrix}
    \bar{x}_{1,f,r} & \bar{x}_{2,f,r}
\end{bmatrix}^T$. We compute from \eqref{lieunicy} that
\begin{align}
    \bar{x}_{i,f,r}(\tau)={}&\bar{x}_{i,f,r}(0)e^{\left(\lambda-\check{k}\alpha\check{\rho}_{x_i}/2\right)\tau}, \quad i=1,2, \\
    \bar{\theta}_{e,r}(\tau)={}&\bar{\theta}_{e,r}(0)-\sum_{i=1}^2\frac{\check{k} \lambda \check{\rho}_{x_i}\bar{x}^2_{i,f,r}(0)}{2(\lambda-\check{k} \alpha\check{\rho}_{x_i}/2)} \nonumber \\
    &\times \left(e^{2\left(\lambda-\check{k} \alpha\check{\rho}_{x_i}/2\right)\tau}-1\right).
\end{align}
Then, we conclude that $(\bar{x}_{1,f,r}, \bar{x}_{2,f,r}, \bar{\theta}_{e,r})=(0, 0, p_{\theta})=:\Upsilon$, where $p_{\theta}=\bar{\theta}_{e,r}(0)+\sum_{i=1}^2\frac{\check{k} \lambda \check{\rho}_{x_i}\bar{x}^2_{i,f,r}(0)}{2(\lambda-\check{k} \alpha\check{\rho}_{x_i}/2)}$, of \eqref{lieunicy} is globally uniformly exponentially stable for $\check{k} \alpha>\frac{2\lambda}{\check{\rho}_{x_i}}$ for $i=1,2$.
This condition implies from \eqref{kcond} and \eqref{rhoxdef} that \eqref{param_unic} needs to be satisfied.

\textbf{Step 4: Singularly perturbed Lie bracket averaging theorem.} 
Combining the global uniform exponential stability of \eqref{lieunicy} at $\Upsilon$ with the global exponential stability of the boundary layer system \eqref{boundlayer} at the origin, by Theorem \ref{thm:singperturb}, we prove the semi-global practical uniform asymptotic stability
of \eqref{xfsingular} at $\Upsilon$. The result holds for $\varepsilon<\varepsilon^*$ with some $\varepsilon^*>0$, which implies $1/\varepsilon=\omega_h-2\lambda>1/\varepsilon^*=:\omega_h^*$. This simplifies to $\omega_h>\omega_h^*+2\lambda$.
The existence of $\omega^*$ and the role of $\omega$ are as defined in Definition \ref{def:singular}.

\textbf{Step 5: Convergence to extremum.} Considering the result in Step 4 and
\begin{align}
    x={}e^{-\lambda t} x_{f}+x^*,
\end{align}
we conclude the semi-global exponential convergence of $x(t)$ to $x^*$ at the rate of $\lambda$. This implies, from \eqref{cost_unic} and \eqref{unicetaf}, the semi-global exponential convergence of the cost function $J(x(t))$ and the filtered state $-\eta(t)$ to $J(x^*)$ at the rate of $2\lambda$. 
\end{proof}

\subsection{Experimental Results}

We demonstrate the effectiveness of the exponential uES designed in \eqref{xdot_unic}--\eqref{thetadot_unic} through a light-seeking experiment for proof-of-concept. 
Fig. \ref{fig:experiment_sub1} depicts the experimental setup, which includes a 750-lumen lamp positioned approximately 80 cm above the surface, and a TurtleBot3 unicycle robot starting about 112 cm from the center of maximum light intensity on the surface. The robot is equipped with a light sensor (Adafruit VEML7700) mounted on top. 
We placed the camera above, parallel to the surface, to capture a top-down view of the robot's motion and used Python's Open Computer Vision (OpenCV) library to visualize the robot's trajectory from the recorded videos.

The parameters in \eqref{xdot_unic}--\eqref{eta_unic} are set to $\omega=0.5, \omega_h=7, k=0.01, \alpha=0.002, \lambda=0.005$, resulting in an initial forward velocity of $u_1(0)=0.1$ m/s and an angular velocity of $u_2(0)=0.5$ rad/s, with $\eta(0)=-J(0)$. The initial position is approximately $x(0)=\begin{bmatrix}
    1.5 & 1.5
\end{bmatrix}^T$ in meters.
For a comparison, we implement classical ES with the same parameters except for $\lambda=0$. 
The vehicle trajectories for both designs are shown in Fig. \ref{fig:experiment_sub1}, where the light intensity shifts from blue to red as it increases. We can see the enhancement of the convergence error with the developed design. The classical ES reaches an average light intensity of 3627 lux at steady state, while  uES increases this value to 4179 lux, exceeding the classical ES by more than 15 percent.
However, due to sensor noise and the large center of maximum light intensity on the surface, further reduction of the convergence error, as theoretically guaranteed, is limited. 

For the classical ES design, reducing the size of the neighborhood around the optimum at steady state is possible by decreasing $\alpha$ in \eqref{xdot_unic} and increasing $k$ in \eqref{thetadot_unic} for the same convergence rate. However, a large $k$ in \eqref{thetadot_unic} can result in a high initial angular velocity, exceeding the robot's speed limits and causing a large deviation from the current path. Our design addresses this issue with smooth transient behavior. A similar discussion is provided in Section \ref{asymuES_numsim1} for Fig. \ref{fig:asymp_uES}.

\section{Conclusion} \label{sec:conclusion}
There has been a long-standing interest in designing ES that can seek and track optima as closely and quickly as possible. We provide various ES designs that achieve unbiased seeking of fixed optima and perfect tracking of time-varying optima at any user-defined rate, including asymptotic, exponential, and in prescribed time. The optima can be static, decaying, or growing unbounded at any rate, including in finite time. These designs leverage methods known as state scaling and time scaling. In essence, these methods use the boundedness of the state-scaled state in a dilated time domain to guarantee convergence of the original system in the original time domain, resulting in accelerated and unbiased optimization. Numerical results are accompanied by experimental results on the problem of source seeking by a unicycle. Future work involves extending these results to achieve unbiased Nash equilibrium seeking.

\section{Acknowledgment}
The authors would like to thank Wenhan Tang, Max Lee, Eric Foss, and Sankalp Kaushik for their tremendous effort in conducting the experiments.

\section*{References}
\vspace{-0.5cm}

\appendix
\subsection{Additional Theorem}
\begin{theorem} \label{generalproof}
Consider the following uES design
\begin{equation} 
\begin{cases}
\begin{aligned}
\dot{\theta}={}&\phi^p(t) \sum_{i=1}^n \sqrt{\alpha_i \omega_i}  e_i  \\
&\times \cos\Big(\omega_i (t_0+\gamma(\phi^q(t)-1))+k_i \phi^r(t)(J(\theta)-\eta)\Big),   \\
\dot{\eta}={}&\left(-\omega_h \eta+\omega_h J(\theta)\right)\phi^{p+1}(t),
\end{aligned}
\end{cases}
\end{equation}
where the function $\phi(t)$ is defined in Table \ref{table_chirpy} with the corresponding parameters $p, q, r, \gamma$. Let $\omega_i = \omega \hat{\omega}_i$ such that $\hat{\omega}_i \neq \hat{\omega}_j$ $\forall i \neq j$, $ t_0 \geq 0$, $\alpha_i, k_i, \omega_h>0$ $\forall i=1,\dots,n$.  
Let Assumptions \ref{Ass0}, \ref{Ass1} hold, and Assumption \ref{asympextbound} holds with $c<p-2\kappa+r$, $d<p-2\kappa+1$. Then, there exists $\omega^*>0$
such that for all $\omega > \omega^*$, the input $\theta(t) \to \theta^*(t)$ semi-globally with respect to $\omega$ and at the same rate that $1/\phi(t) \to 0$, and there exist a class $\mathcal{KL}$ function $\mathcal{B}$, a class $\mathcal{K}$ function $\mathcal{Y}$, and a nonnegative constant $D(\theta(t_0), \theta^*(t_0), \eta(t_0))$ such that 
    \begin{align}
    |\theta(t)-\theta^*(t)|\leq{}&\phi^{-1}(t) \Big(D+\mathcal{B}\left(|\theta(t_0)-\theta^*(t_0)|, t-t_0\right) \nonumber \\
    &+\sup_{t_0 \leq s \leq t}\mathcal{B}\left(\mathcal{Y}\left( \phi^{-1}(s) \right), t-s\right)\Big). \label{KLforgeneral}
\end{align}
\end{theorem}

\begin{proof}
\textbf{Step 1: State transformation.}
Let us consider the following transformations
\begin{align}
    \theta_{f}={}& \phi(t) (\theta-\theta^*(t)), \label{statetransfor_app} \\
    \eta_{f}={}& \phi^{2\kappa}(t) (\eta-J(\theta^*(t))), \label{statetransfor2_app}
\end{align}
which transform \eqref{ESasymp}  to
\begin{equation} \label{thetaasymtransappendix}
\begin{cases}
\begin{aligned}
\dot{\theta}_{f}={}&-\phi(t)\dot{\theta}^*(t)+\frac{\dot{\phi}(t)}{\phi(t)}\theta_f+\phi^{p+1}(t) \sum_{i=1}^n \sqrt{\alpha_i \omega_i}  e_i   \\
&\times \cos\Big(\omega_i (t_0+\gamma(\phi^q(t)-1))+k_i \phi^r(t) \\
&\times \left(J_f(\theta_f, t)-\phi^{-2\kappa}(t){\eta}_f\right)\Big),  \qquad \quad \\
\dot{\eta}_f={}&\left(\frac{2\kappa  \dot{\phi}(t) }{\phi(t)}-\omega_h \phi^{p+1}(t) \right)\eta_f+\omega_h \phi^{2\kappa+p+1}(t) J_f(\theta_f, t)\\
&-\phi^{2\kappa}(t)\dot{J}(\theta^*(t)), 
 \end{aligned}
\end{cases}
\end{equation}
with
\begin{align}  
    J_f(\theta_f, t) = {}& J(\theta_f/\phi(t)+\theta^*(t)) -J(\theta^*(t)).
\end{align}
\textbf{Step 2: Time transformation.} We introduce the following time dilation and contraction transformations
\begin{align}
    \tau={}&t_0+\gamma(\phi^q(t)-1) \label{taudefapp}\\
    t={}&\phi^{-1}\Big(\left(1+(\tau-t_0)/\gamma\right)^{\frac{1}{q}}\Big), \label{timecontractgen}
\end{align}
where $\phi^{-1}(\cdot)$ denotes the inverse of the function $\phi(\cdot)$. Each $\phi(t)$ function in Table \ref{table_chirpy} with the corresponding the parameters $p, q$ yields
\begin{align}
    \dot{\phi}(t)={}&\phi^{p-q+2}(t)/(\gamma q). \label{phidotapp}
\end{align}
Using \eqref{taudefapp} and \eqref{phidotapp}, we compute
\begin{align}
    \frac{d\tau}{dt}={}&\phi^{p+1}(t).
\end{align}
Then, we rewrite \eqref{thetaasymtransappendix} in dilated $\tau$-domain as
\begin{equation} \label{dilatedthetafapp}
\begin{cases}
\begin{aligned}
    \frac{d\theta_f}{d\tau}={}&-\phi_{\tau}(\tau) \frac{d\theta_{\tau}^*(\tau)}{d\tau}+\frac{1}{\gamma q \phi^q_{\tau}(\tau)}\theta_f+\sum_{i=1}^n \sqrt{\alpha_i \omega_i}e_i  \\
    &\times \cos\Big(\omega_i \tau +k_i \phi_{\tau}^r(\tau)(J_{f,\tau}(\theta_f, \tau)-\phi^{-2\kappa}_{\tau}(\tau)\eta_f)\Big) \\
    \frac{d\eta_f}{d\tau}={}&\left(\frac{2\kappa}{\gamma q \phi^q_{\tau}(\tau)}-\omega_h\right) \eta_f+\omega_h \phi_{\tau}^{2\kappa}(\tau) J_{f,\tau}({\theta}_f, \tau)\\
    &-\phi_{\tau}^{2\kappa}(\tau)\frac{\partial J(\theta_{\tau}^*(\tau))}{\partial \tau},
\end{aligned}
\end{cases}
\end{equation}
where
\begin{align}
    \phi_{\tau}(\tau)={}&(1+(\tau-t_0)/\gamma)^{\frac{1}{q}}, \label{phitaudef} \\
    \theta_{\tau}^*(\tau)={}&\theta^*\left(\phi^{-1}\left(\left(1+(\tau-t_0)/\gamma\right)^{\frac{1}{q}}\right)\right)
\end{align}
with
\begin{align}
    J_{f,\tau}(\theta_f,\tau)={}J(\theta_f/\phi_{\tau}(\tau)+\theta_{\tau}^*(\tau)) -J(\theta_{\tau}^*(\tau)).
\end{align}
\textbf{Step 4: Stability analysis.} Note that the system \eqref{dilatedthetafapp} is in similar form to \eqref{thetaasymtrans} with $\gamma=\frac{1}{\beta}, q=v$, except that the bounds in Assumption \ref{asympextbound} need to be satisfied in $\tau$-domain for \eqref{dilatedthetafapp}. The growth bounds in Assumption \ref{asympextbound} is rewritten in $\tau$-domain as
\begin{align}
\left|\frac{d{\theta}_{\tau}^*(\tau)}{d\tau}\right|+\phi^{p+1}_{\tau}(\tau)\left|\frac{d^2{\theta}^*(t)}{d\tau^2}\right|\leq {} & M_{\theta} \phi_{\tau}^{c-p-1}(\tau),  \\
\left| \frac{\partial J(\theta_{\tau}^*(\tau))}{\partial\tau}\right|+\phi^{p+1}_{\tau}(\tau)\left|\frac{\partial^2 J(\theta_{\tau}^*(\tau))}{\partial\tau^2} \right|\leq {}&  M_J\phi_{\tau}^{d-p-1}(\tau),
\end{align}
for $\tau \in [t_0, \infty)$.
Recall that the conditions on the powers $c, d$ in Assumption \ref{asympextbound}, for the practical stability of the system \eqref{thetaasymtrans} is that $c<-1-2\kappa+r$ and $d<-2\kappa$, as provided in Table \ref{table_consfreq}. Following this analogy, for the stability of \eqref{dilatedthetafapp}, the powers of the function $\phi_{\tau}(\tau)$ should satisfy
\begin{align}
    c-p-1<&-1-2\kappa+r, \\
    d-p-1<&-2\kappa,
\end{align}
which are simply 
\begin{align}
    c<{}&p-2\kappa+r, \\
    d<{}&p-2\kappa+1.
\end{align}
Following the Step 2 to 6 in the proof of Theorem \ref{theoremasymp}, we prove that the origin of \eqref{dilatedthetafapp} is practically uniformly asymptotically stable. 
By the time contraction \eqref{timecontractgen}, we conclude the practical uniform stability of \eqref{thetaasymtransappendix}
with the existence of $\omega^*$ and the role of $\omega$ attributed as in Definition \ref{def:localpractical}.

Similar to \eqref{ISSthetaf}, the fading memory ISS bound of the average state $\bar{\theta}_f$ in $t$-domain is obtained as
\begin{align}
    |\bar{\theta}_f(t)| \leq \mathcal{B}\left(|\bar{\theta}_f(t_0)|, t-t_0\right)+\sup_{t_0 \leq s \leq t}\mathcal{B}\left(\mathcal{Y}\left( \phi^{-1}(s) \right), t-s\right),
\end{align}
where the functions $\mathcal{B}$ and $\mathcal{Y}$ are defined as in \eqref{ISSbound1asympt} and \eqref{ISSbound2asympt}, except that $\xi(t), \beta$, and $v$ in \eqref{ISSbound1asympt} and \eqref{ISSbound2asympt} are replaced by $\phi(t), \frac{1}{\gamma}$, and $q$, respectively.

\textbf{Step 5: Convergence to extremum.} Considering the result in Step 4 and recalling from \eqref{statetransfor_app} that
\begin{align}
    \theta=\theta^*(t)+\frac{1}{\phi(t)} \theta_f, 
\end{align}
we conclude the convergence of $\theta(t)$ to $\theta^*(t)$ at the same rate that $1/\phi(t) \to 0$. 
The bound on the convergence error is obtained  as 
\begin{align}
    |\theta(t)-\theta^*(t)| \leq{} \phi^{-1}(t)\left(|\bar{\theta}_f(t)|+|\theta_f(t)-\bar{\theta}_f(t)|\right), 
\end{align}
from which we conclude \eqref{KLforgeneral}.
From \eqref{Jboundcon}, \eqref{statetransfor_app}, and \eqref{statetransfor2_app}, we prove the convergence of the output $y(t)$ and the filtered state $\eta(t)$ to $J(\theta^*(t))$ at the rate of $1/\phi^{2\kappa}(t)$
and complete the proof of Theorem \ref{generalproof}.
\end{proof}

\subsection{Useful Definitions} \label{sec:app2}
\begin{definition}[\cite{durr2013lie}] \label{def:localpractical}
A compact set $\mathcal{S} \subset \mathbb{R}^n $ is said to be \textbf{locally practically uniformly asymptotically stable} for \eqref{conaff} if the following three conditions are satisfied:
\begin{itemize}
    \item \underline{Practical Uniform Stability}: For any $\epsilon>0$ there exist
$\delta, \omega^*>0$ such that for all $t_0 \in \mathbb{R}^+$ and $\omega> \omega^*$, if
    $x(t_0) \in U_{\delta}^{\mathcal{S}}$, then $x(t) \in U_{\epsilon}^{\mathcal{S}}$ for $t \in [t_0, \infty)$.
    \item \underline{$\delta$-Practical Uniform Attractivity}:  Let $\delta>0$. For any $\epsilon>0$ there exist $t_1 \geq 0$ and $\omega^* > 0$ such that for all $t_0 \in \mathbb{R}^+$ and $\omega>\omega^*$, if
    $x(t_0) \in U_{\delta}^{\mathcal{S}}$, then $x(t) \in U_{\epsilon}^{\mathcal{S}}$ for $ t \in [t_0+t_1, \infty)$.
    \item \underline{Practical Uniform Boundedness}:  For any $\delta>0$ there exist $\epsilon>0$ and $\omega^*>0$ such that for all $t_0 \in \mathbb{R}^+$ and $\omega>\omega^*$, if
    $x(t_0) \in U_{\delta}^{\mathcal{S}}$, then $x(t) \in U_{\epsilon}^{\mathcal{S}}$ for $t \in [t_0, \infty)$.
\end{itemize}
Furthermore, if $\delta$-practical uniform attractivity holds for every $\delta>0$, then the compact set $\mathcal{S}$ is said to be \textbf{semi-globally practically uniformly asymptotically stable} for \eqref{conaff}.

\end{definition}

\begin{definition}[\cite{durr2017extremum}] \label{def:singular}
A compact set $\mathcal{S} \subset \mathbb{R}^n $ is said to be \textbf{singularly semi-globally practically uniformly asymptotically stable} for \eqref{singpert1} if the following three conditions are satisfied:
\begin{itemize}
    \item \underline{Singular Practical Uniform Stability}: For all $\epsilon_x, \epsilon_z \in (0, \infty)$ there exist $\delta_x, \delta_z \in (0, \infty)$ and $\omega^* \in (0, \infty)$ such that for all $\omega \in (\omega^*, \infty)$ there exists a $\varepsilon^* \in (0, \infty)$ such that for all $\varepsilon \in (0, \varepsilon^*)$ and for all $t_0 \in \mathbb{R}$ 
if $x_0 \in U_{\delta_x}^{\mathcal{S}}$ and $z_0 - l(x_0) \in U^0_{\delta_z} $, then
    $x(t;t_0, x_0) \in U^S_{\epsilon_x}$ and $z(t;t_0, z_0) - l(x(t;t_0, x_0)) \in U^0_{\epsilon_z}$ for $t \in [t_0,\infty)$.
    \item \underline{Singular Practical Uniform Attractivity}: For all $\delta_x$, $\delta_z$, $\epsilon_x$, $\epsilon_z \in (0, \infty)$, there exist a $t_f \in [0,\infty)$ and $\omega^* \in (0,\infty)$ such that for all $\omega \in (\omega^*,\infty)$ there exists a $\varepsilon^* \in (0,\infty)$ such that for all $\varepsilon \in (0, \varepsilon^*)$ and for all $t_0 \in \mathbb{R}$, if
    $x_0 \in U_{\delta_x}^{\mathcal{S}}$ and $z_0 - l(x_0) \in U^0_{\delta_z}$, then $x(t;t_0, x_0) \in U^S_{\epsilon_x}$, $t \in [t_0+\frac{t_f}{\varepsilon},\infty)$ and $z(t;t_0, z_0)- l(x(t;t_0, x_0)) \in U^0_{\epsilon_z}$ for $t \in [t_0+t_f,\infty)$.
    \item \underline{Singular Practical Uniform Boundedness}: For all $\delta_x$, $\delta_z \in (0,\infty)$ there exist $\epsilon_x, \epsilon_z \in (0,\infty)$ and $\omega^* \in (0,\infty)$ such that for all $\omega \in (\omega^*,\infty)$ there exists a $\varepsilon^* \in (0,\infty)$ such that for all $\varepsilon \in (0, \varepsilon^*)$ and for all $t_0 \in \mathbb{R}$, if $x_0 \in U_{\delta_x}^{\mathcal{S}}$ and $z_0-l(x_0) \in U^0_{\delta_z}$, then $x(t;t_0, x_0) \in U^S_{\epsilon_x}$ and $z(t;t_0, z_0) - l(x(t;t_0, x_0)) \in U^0_{\epsilon_z}$ for  $t \in [t_0,\infty)$.
\end{itemize}
\end{definition}

\end{document}